\documentclass[12pt,oneside]{article}
 \newtheorem{thm}{Theorem}[section]
 \newtheorem{cor}[thm]{Corollary}
 
 \newtheorem{propo}[thm]{Proposition}
 
\usepackage[T1]{fontenc}
\usepackage{amsmath}
\usepackage{amssymb}
\usepackage{tabularx}
\usepackage[dvips]{epsfig}
\usepackage[dvips]{graphics}
\usepackage{graphicx}
\usepackage{dsfont}

\newcommand{\Prob}{\mathds{P}}
\newcommand{\Nor}{\mathcal{N}}

\newcommand{\Inte}{\mathds{N}}
\newcommand{\Inter}{\mathds{Z}}
\newcommand{\LL}{\mathds{L}}
\newcommand{\Real}{\mathds{R}}
\newcommand{\Esp}{\mathds{E}}
\newcommand{\Var}{\mbox{Var}}

\newcommand{\Proba}{\mathcal{P}}
\newcommand{\Loi}{\mathcal{L}}

\newtheorem{theorem}{Theorem}[section]

\newtheorem{lemma}[theorem]{Lemma}

\newtheorem{proposition}[theorem]{Proposition}

\newtheorem{dem*}[theorem]{Proof}

\def\limitedistl{\renewcommand{\arraystretch}{0.5}
\begin{array}[t]{c}
\stackrel{{\Loi}}{\longrightarrow} \\
{\scriptstyle N\rightarrow\infty}
\end{array}\renewcommand{\arraystretch}{1}}

\def\limiteprob{\renewcommand{\arraystretch}{0.5}
\begin{array}[t]{c}
\stackrel{{\Proba}}{\longrightarrow} \\
{\scriptstyle N\rightarrow\infty}
\end{array}\renewcommand{\arraystretch}{1}}

\def\limiteNN{\renewcommand{\arraystretch}{0.5}
\begin{array}[t]{c}
\stackrel{}{\longrightarrow} \\
{\scriptstyle N\rightarrow\infty}
\end{array}\renewcommand{\arraystretch}{1}}

\begin{document}
\date{}
\title{Detecting changes in the fluctuations of a Gaussian process and an application
to heartbeat time series} \maketitle \vspace{-1cm}
\begin{center}
\author{Jean-Marc Bardet~~~~Imen Kammoun}\\
~\\
Université Paris 1, SAMOS-MATISSE-CES, 90 rue de Tolbiac, 75013
Paris, France.\\
\texttt{jean-marc.bardet@univ-paris1.fr},~\texttt{imen.kammoun@univ-paris1.fr}\\
~\\
~\\
\end{center}
\begin{abstract}
The aim of this paper is first the detection of multiple abrupt
changes of the long-range dependence (respectively self-similarity,
local fractality) parameters from a sample of a Gaussian stationary
times series (respectively time series, continuous-time process
having stationary increments). The estimator of the $m$ change
instants (the number $m$ is supposed to be known) is proved to
satisfied a limit theorem with an explicit convergence rate.
Moreover, a central limit theorem is established for an estimator of
each long-range dependence (respectively self-similarity, local
fractality) parameter. Finally, a goodness-of-fit test is also built
in each time domain without change and proved to asymptotically
follow a Khi-square distribution. Such statistics are applied to
heart rate data of marathon's runners and lead to interesting
conclusions.
\end{abstract}
{\em Keywords:}  Long-range dependent processes; Self-similar
processes; Detection of abrupt changes; Hurst parameter;
Self-similarity parameter; Wavelet analysis; Goodness-of-fit test.
\section{Introduction}
\begin{center}
\begin{figure}[h]
\begin{center}
\includegraphics[width=9 cm,height=5 cm]{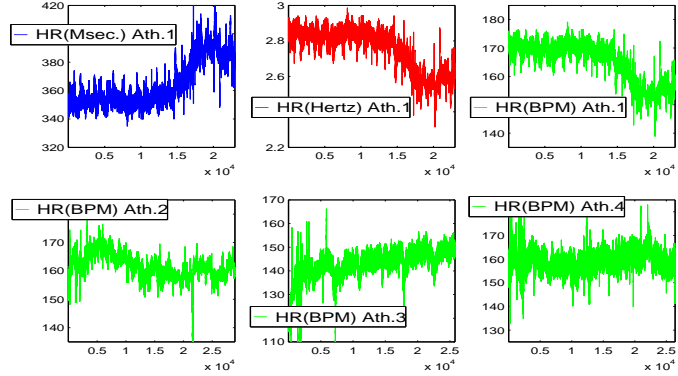}
\end{center}
\caption{Heat rate signals of Athlete 1 in ms, Hertz and BPM (up),
of Athletes 2, 3 and 4 in BPM (down)} \label{Figure1}
\end{figure}
\end{center}

The  content of this paper was motivated by a
general study of physiological signals of runners recorded during
endurance races as marathons. More precisely, after different signal
procedures for "cleaning" data, one considers the time series
resulting of the evolution of heart rate (HR data in the sequel)
during the race. The following figure provides several examples of
such data (recorded during Marathon of Paris 2004 by Professor V.
Billat and her laboratory LEPHE, see http://www.billat.net). For
each runner, the periods (in ms) between the successive pulsations
(see Fig. \ref{Figure1}) are recorded. The HR signal in number of
beats per minute (bpm) is then deduced (the HR average for the whole
sample is of 162 bpm).
Numerous authors have studied heartbeat time series (see for
instance \cite{pengh}, \cite{pengm} or \cite{absil}). A model
proposed to fit these data is a trended long memory process with an
estimated Hurst parameter close to $1$ (and sometimes more than
$1$). In \cite{bardet2} three improvements have been proposed to
such a model: 1/ data are stepped in three different stages which
are detected using a change point's detection method (see for
instance \cite{lavielle1} or \cite{lavielle2}). The main idea of the
detection's method is to consider that the signal distribution
depends on a vector of unknown characteristic parameters constituted
by the mean and the variance. The different stages (beginning,
middle and end of the race) and therefore the different vectors of
parameters, which change at two unknown instants, are estimated. 2/
during each stage, a time-continuous Gaussian process is proposed
for modelling the detrended time series. This process is a
generalization of a fractional Gaussian noise (FGN) also called
locally fractional Gaussian noise such that, roughly speaking, there
exists a local-fractality parameter $H\in \Real$ (corresponding to
Hurst parameter for FGN) only for frequencies $|\xi| \in
[f_{min}\,,\, f_{max}]$ with $0<f_{min}<f_{max}$ (see more details
below). 3/ this parameter $H$ which is very interesting for
interpreting and explaining the physiological signal behaviours, is
estimating from a wavelet analysis. Rigorous results are also proved
providing a
central limit theorem satisfied by the estimator. \\
In order to improve this study of HR data and since the eventual
changes of $H$ values are extremely meaningful for explaining the
eventual physiological changes of the athlete's HR during the race,
the detection of abrupt change of $H$ values is the aim of this
paper. By this way the different stages detected during the race
will be more relevant for explaining the physiological status of the
athlete than stages detected from changes in mean or variance. For
instance, the HR of a runner could decrease in mean
even if the "fluctuations" of the HR does not change.\\
~\\
In this paper, an estimator of $m$ instants ($m \in \Inte^*$) of
abrupt changes of long-range dependence, self-similarity or
local-fractality (more details about these terms will be provided
below) is developed for a sample of a Gaussian process. Roughly
speaking, the principle of such estimator is the following: in each
time's domain without change, the parameter of long-range dependence
(or self-similarity or local self-fractality) can be estimated from
a log-log regression of wavelet coefficients' variance onto several
chosen scales. Then a contrast defined by the sum on every $m+1$
possible zones of square distances between points and regressions
lines is minimized providing an estimator of the $m$ instants of
change. Under general assumptions, a limit theorem with a
convergence rate satisfied by such an estimator is
established in Theorem \ref{convprob}. \\
Moreover, in each estimated no-change zone, parameters of long-range
dependence (or self-similarity or local self-similarity) can be
estimated, first with an ordinary least square (OLS) regression,
secondly with a feasible generalized least square (FGLS) regression.
Central limit theorems are established for both these estimators
(see Theorem \ref{convparam} and Proposition \ref{FGLS} below) and
confidence intervals can therefore be computed. The FGLS estimator
provides two advantages: from the one hand, its asymptotic variance
is smaller than OLS estimator one. From the other hand, it allows to
construct a very simple (Khi-square) goodness-of-fit test based on a
squared distance between points and FGLS regression line. The
asymptotic
behavior of this test is provided in Theorem \ref{T}. \\
~\\
Then, different particular cases of Gaussian processes are studied:
\begin{enumerate}
\item long-range dependent processes with abrupt changes of
values of LRD parameters. In such time series case, a
semi-parametric frame is supposed (including fractional Gaussian
noises (FGN) and Gaussian FARIMA processes) and assumptions of limit
theorems are always satisfied with interesting convergence rates
(see Corollary \ref{cor_LRD}).
\item self-similar time series with abrupt changes of
values of self-similarity parameters. In such case, fractional
Brownian motions (FBM) are only considered. Surprisingly,
convergences of estimators are only established when the maximum of
differences between self-similarity parameters is sufficiently
small. Simulations exhibit a non convergence of the estimator of
instant change when a difference between two parameters is too large
(see Corollary \ref{cor_FBM}).
\item locally fractional Gaussian processes with abrupt changes
of values of local-fractality parameters. In such a continuous time
processes' case, a semi-parametric frame is supposed (including
multiscale fractional Brownian motions) and assumptions of limit
theorems are always satisfied with interesting convergence rates
(see Corollary \ref{cor_local}).
\end{enumerate}
The problem of change-point detection using a contrast minimization
was first studied in the case of independent processes (see for
instance Bai and Perron \cite{Bai-Per}), then for weakly dependent
processes (see for instance Bai \cite{bai}, Lavielle
\cite{lavielle1} or Lavielle and Moulines \cite{Lav-Mou}) and since
middle of 90's in the case of processes which exhibit long-range
dependance (see for instance Giraitis {\it et al.} \cite{gls},
Kokoszka and Leipus \cite{kl} or Lavielle and Teyssi\`ere
\cite{lavielle2}). Of the various approaches, some were associated
with a parametric framework for a change points detection in mean
and/or variance and others where associated with a non-parametric
framework (typically like detecting changes in distribution or
spectrum). To our knowledge, the semi-parametric case of abrupt
change detection for long-range dependent or self-similarity
parameter
is treated here for the first time.  \\
However, in the literature different authors have proposed test
statistics for testing the no-change null hypothesis against the
alternative that the long-memory parameter changes somewhere in the
observed time series. Beran and Terrin \cite{bt} proposed an
approach based on the Whittle estimator, Horv\'ath and Shao
\cite{hs} obtained limit distribution of the test statistic based on
quadratic forms and Horv\'ath \cite{hor} suggested another test
based on quadratic forms of Whittle estimator of long-memory
parameter. The goodness-of-fit test presented below and which
satisfies the limit theorem \ref{T} also allows to test if the
long-range memory (or self-similarity or local-fractality)
parameter changes somewhere in the time series. \\
~\\
Our approach is based on the wavelet analysis. This method applied
to LRD or self-similar processes for respectively estimating the
Hurst or self-similarity parameter was introduced by Flandrin
\cite{flan2} and was developed by Abry, Veitch and Flandrin
\cite{avf} and Bardet {\it et al.} \cite{JM2}. The convergence of
wavelet analysis estimator was studied in the case of a sample of
FBM in \cite{JM1}, and in a semi-parametric frame of a general class
of stationary Gaussian LRD processes by Moulines {\it et al.}
\cite{mrt} and  Bardet {\it et al.} \cite{JM2}. Moreover, wavelet
based estimators are robust in case of polynomial trended processes
(see Corollary \ref{trend}) and is therefore very interesting for
studying stochastic fluctuations of a process
without taking care on its smooth variations.\\
A method based on wavelet analysis was also developed by Bardet and
Bertrand \cite{bardet1} in the case of multiscale FBM (a
generalization of the FBM for which the Hurst parameter depends on
the frequency as a piecewise constant function) providing statistics
for the identification (estimation and goodness-of-fit test) of such
a process. Such a process was used for modelling biomechanics
signals. In the same way, the locally fractional Gaussian process (a
generalization of the FBM for which the Hurst parameter, called the
local-fractality parameter, is constant in a given domain of
frequencies) was studied in \cite{bardet2} for modelling HR data
during the three characteristics stages of the race. An increasing
evolution of the local-fractality parameter during the race was
generally showed for any runner from this method. Using the method
of abrupt change detection of local-fractality parameter $H$
developed in Corollary \ref{cor_local}, this result is confirmed by
estimations of $H$ for each runner even if the change's instants
seem to vary a lot depending on the fatigue of the runner (see the
application to
HR's time series in Section \ref{applications}).\\
~\\
The paper is organized as follows. In Section \ref{results},
notations, assumptions and limit theorems are provided in a general
frame. In Section \ref{applications}, applications of the limit
theorems to three kind of "piecewise" Gaussian process are presented
with also simulations. The case of HR data is also treated. Section
\ref{proofs} is devoted to the proofs.
\section{Main results}\label{results}
\subsection{Notations and assumptions}
First, a general and formal frame can be proposed. Let $(X_t)_{t \in
T}$ be a zero-mean Gaussian process with $T=\Inte$ or $T=\Real$ and
assume that
$$
\big (X_0,X_{\delta_N},X_{2\delta_N},\ldots,X_{N \delta_N}\big
)~~\mbox{is known with $\delta_N=1$ or $\delta_N \limiteNN 0$,}
$$
following data are modeled with a time series ($T=\Inte$) or a
continuous time process $T=\Real$. In the different proposed
examples $X$ could be a stationary long memory time series or a
self-similar or locally fractional process having stationary increments. \\
For estimations using a wavelet based analysis, consider $\psi:\Real
\rightarrow \Real$ a function called "the mother wavelet". In
applications, $\psi$ is a function with a compact (for instance
Daubeshies wavelets) or an essentially compact support (for instance
Lemari\'e-Meyer wavelets). For $(X_t)_{t \in T}$ and $(a,b)\in
\Real_+^*\times \Real$, the wavelet coefficient of $X$ for the scale
$a$ and the shift $b$ is $$
d_{X}(a,b):=\frac{1}{\sqrt{a}}\int_{\Real}\psi(\frac{t-b}{a})X(t)dt.$$
When only a discretized path of $X$ is available (or when $T
=\Inte$), approximations $e_{X}(a,b)$ of $d_{X}(a,b)$ are only
computable. We have chosen to consider for $(a,b)\in \Real_+^*\times
\Inte$,
\begin{eqnarray}\label{coeff_e}
e_{X}(a,b):=\frac {\delta_n}{\sqrt a} \, \sum_{p=1}^{N} \psi\big
(\frac {p-b} {a} \big ) X_{p \, \delta_N  },
\end{eqnarray}
which is the formula of wavelet coefficients computed from Mallat's
algorithm for compactly supported discrete ($a\in 2^\Inte$) wavelet
transform (for instance Daubeshies wavelets) when $N$ is large
enough and nearly this formula for discrete wavelet transform with
an essentially compact support (for instance Lemari\'e-Meyer
wavelets). Now assume that there exist $m\in \Inte$ (the number of
abrupt changes) and
\begin{itemize}
\item $0=\tau_0^*<\tau_1^*<\ldots<\tau_{m}^*<\tau_{m+1}^*=1$ (unknown parameters);
\item two families $(\alpha_j^*)_{0\leq j \leq m}\in \Real^{m+1}$ and
$(\beta_j^*)_{0\leq j \leq m}\in (0,\infty)^{m+1}$ (unknown
parameters);
\item a sequence of "scales" $(a_n)_{n \in \Inte}\in \Real^\Inte$ (known sequence) satisfying
$a_n \geq a_{min}$ for all $n \in \Inte$, with $a_{min}>0$,
\end{itemize}
such that for all $j=0,1,\ldots,m$ and $k \in D_N^*(j) \subset \big
[[N\delta_N\tau_j^*]\, , \,[N\delta_N\tau_{j+1}^*]\big ]$,
\begin{eqnarray}\label{esper_e}
\Esp \big [e_{X}^2(a_N,k)\big ] \sim \beta_j^* \cdot \big ( a_N\big
)^{\alpha_j^*}~~\mbox{when}~~N\to \infty \mbox{ and } N\delta_N\to
\infty.
\end{eqnarray}
Roughly speaking, for $N\in \Inte^*$ the change instants are
$[N\delta_N\tau_j^*]$ for $j=1,\ldots,m$, the variance of wavelet
coefficients follows a power law of the scale, and this power law is
piecewise varying following the shift. Thus piecewise sample
variances can be appropriated estimators of parameters of these
power laws. Hence let us define
\begin{eqnarray}\label{Sk}
S^{k'}_{k}(a_N):= \frac {a_N} {k'-k} \, \sum
_{p=[k/a_N]}^{[k'/a_N]-1}e_{X}^2(a_N,a_N\,p)~~\mbox{for $0\leq k <
k'\leq N\delta_N$}.
\end{eqnarray}
Now set $0<r_1<\ldots<r_\ell$ with $\ell \in \Inte^*$ and let us
suppose that a multidimensional central limit theorem can also be
established for $\big(S^{k'}_{k}(r_i \, a_N)\big )_{1 \leq i \leq
\ell}$, {\it i.e.}
\begin{eqnarray}\label{TLCS}
\big(S^{k'}_{k}(r_i \, a_N)\big )_{1 \leq i \leq \ell}=\big (
\beta_j^* \cdot \big (r_i \,  a_N\big )^{\alpha_j^*}\big )_{1 \leq i
\leq \ell} +\big (a_N\big )^{\alpha_j^*} \times \sqrt {\frac {a_N
}{k'-k}} \, \big ( \varepsilon_i^{(N)}(k,k') \big )_{1 \leq i \leq
\ell},
\end{eqnarray}
with $[N\delta_N\tau_j^*] \leq k < k' \leq [N\delta_N\tau_{j+1}^*]$
and it exists $\Gamma^{(j)}(\alpha_j^*,r_1,\ldots,r_\ell)=\big (
\gamma_{pq}^{(j)}\big ) _{1 \leq p,q \leq \ell}$ a $(\ell \times
\ell)$ matrix not depending on $N$ such that $\alpha \mapsto
\Gamma^{(j)}(\alpha,r_1,\ldots,r_\ell)$ is a continuous function, a
positive matrix for all $\alpha$ and
\begin{eqnarray}\label{lim_epsi}
\big ( \varepsilon_i^{(N)}(k,k') \big )_{1 \leq i \leq \ell}
\limitedistl
\Nor\big(0,\Gamma^{(j)}(\alpha_j^*,r_1,\ldots,r_\ell)\big)~~\mbox{when}~~k'-k
\to \infty.
\end{eqnarray}
With the usual Delta-Method, relation (\ref{TLCS}) implies that for
$1 \leq i \leq \ell$,
\begin{eqnarray}\label{TLClogS}
\log \big (S^{k'}_{k}(r_i \, a_N)\big )=\log (\beta_j^*)+\alpha_j^*
\log \big (r_i \, a_N\big )+\sqrt {\frac {a_N }{k'-k}} \,
\varepsilon_i^{(N)}(k,k'),
\end{eqnarray}
for $[N\delta_N\tau_j^*] \leq k < k' \leq [N\delta_N\tau_{j+1}^*]$
and the limit theorem (\ref{lim_epsi}) also holds. This is a linear
model and therefore a log-log regression of $\big (S^{k'}_{k}(r_i \,
a_N)\big )_i$ onto
$\big (r_i \, a_N\big)_i$ provides an estimator of $\alpha_j^*$ and $\log (\beta_j^*)$.\\
~\\
The first aim of this paper is the estimation of unknown parameters
$(\tau_j^*)_j$, $(\alpha_j^*)_j$ and $(\beta_j^*)_j$. Therefore,
define a contrast function
\begin{multline*}
U_N\big((\alpha_j)_{0\leq j\leq m},\,(\beta_j)_{0\leq j\leq m},\,
(k_j)_{1\leq j\leq m}\big)=\sum_{j=0}^{m}\sum_{i=1}^{\ell} \Big(
\log\big(S^{k_{j+1}}_{k_j}(r_i \,a_N)\big)- \big(\alpha_j\log (r_i
\,a_N)+\log \beta_j\big )\Big)^2 \\
\mbox{with} ~~~~ \left \{
\begin{array}{l} \bullet ~
(\alpha_j)_{0\leq j\leq m}\in A^{m+1}\subset \Real^{m+1}\\
\bullet ~(\beta_j)_{0\leq j\leq m} \in B^{m+1} \subset
(0,\infty)^{m+1}\\
\bullet ~0=k_0<k_1<\ldots<k_m<k_{m+1}=N\delta_N, (k_j)_{1\leq j\leq
m} \in K_{m}(N) \subset \Real^{m}
\end{array}\right . .
\end{multline*}
The vector of estimated parameters $\widehat \alpha_j,~\widehat
\beta_j$ and $\widehat k_j$ (and therefore $\widehat \tau_j$) is the
vector which minimizes this contrast function, {\it i.e.},
\begin{eqnarray}
\nonumber \big ((\widehat\alpha_j)_{0\leq j\leq m},\,(\widehat
\beta_j)_{0\leq j\leq m},\, (\widehat k_j)_{1\leq j\leq m}\big
)&& \\
&&\hspace{-6cm}:=\text{Argmin} \Big\{U_N\big((\alpha_j)_{0\leq j\leq
m},\,(\beta_j)_{0\leq j\leq m},\, (k_j)_{1\leq j\leq
m}\big)\Big\}~~\mbox{in}~~A^{m+1}\times B^{m+1}\times K_{m}(N)\\
\hspace{-5.5cm}\widehat \tau_j\hspace{5.5cm}&&
\hspace{-6cm}:=\widehat k_j/(N\delta_N)~~\mbox{for}~~1\leq j\leq m.
\end{eqnarray}
For a given $(k_j)_{1\leq j\leq m}$, it is obvious that
$(\widehat\alpha_j)_{0\leq j\leq m}$ and $( \log
\widehat\beta_j)_{0\leq j\leq m}$ are obtained from a log-log
regression of $\big (S^{k_{j+1}}_{k_j}(r_i \, a_N)\big )_i$ onto
$\big (r_i \, a_N\big)_i$, {\it i.e.}
$$
\left ( \begin{array}{c} \widehat \alpha_j \\
\log \widehat \beta_j \end{array} \right) =\big ( L_{1}' \cdot
L_{1})^{-1} L_{1}'\cdot Y_{k_j}^{k_{j+1}}
$$
with
$\displaystyle{Y_{k_j}^{k_{j+1}}:=\big(\log\big(S^{k_{j+1}}_{k_j}(r_i\cdot
a_N)\big)\big)_{1\leq i\leq \ell}}$ and
$\displaystyle{L_{a_N}:=\left (
\begin{array}{cc} \log (r_1\, a_N) & 1\\ \vdots & \vdots \\\log (r_\ell \, a_N) & 1
\end{array} \right).}$
Therefore the estimator of the vector $(k_j)_{1\leq j\leq m}$ is
obtained from the minimization of the contrast
\begin{eqnarray}
G_N(k_1,k_2,\ldots,k_m)&:=&U_N\big((\widehat \alpha_j)_{0\leq j\leq
m},\,(\widehat \beta_j)_{0\leq j\leq m},\, (k_j)_{1\leq
j\leq m}\big)\\
&& \hspace{-3cm}\Longrightarrow~~(\widehat k_j)_{1\leq j\leq
m}=\text{Argmin} \Big\{G_N(k_1,k_2,\ldots,k_m),~~(k_j)_{1\leq j\leq
m} \in K_{m}(N)\Big \}.
\end{eqnarray}

\subsection{Estimation of abrupt change time-instants $(\tau_j^*)_{1\leq j \leq m}$}
In this paper, parameters $(\alpha_j^*)$ are supposed to satisfied
abrupt changes. Such an hypothesis is provided by the following
assumption:\\
~\\
{\bf Assumption C:} Parameters $(\alpha_j^*)$ are such that $
|\alpha_{j+1}^*-\alpha_j^*|\neq 0~~\mbox{for
all}~~j=0,1,\ldots,m-1.$\\
~\\
Now let us define:
$$
\underline{\tau}^*:=(\tau^*_1,\ldots,\tau^*_m),~~\underline{\widehat\tau}:=(\widehat
\tau_1,\ldots,\widehat
\tau_m)~~\mbox{and}~~\|\underline{\tau}\|_m:=\max \big
(|\tau_1|,\ldots,|\tau_m|\big ).
$$
Then $\underline{\widehat \tau}$ converges in probability to
$\underline{\tau}^*$ and more precisely,
\begin{theorem}\label{convprob} Let $\ell \in \Inte \setminus \{0,1,2\}$.
If Assumption C and relations (\ref{TLCS}), (\ref{lim_epsi}) and
(\ref{TLClogS}) hold with $(\alpha_j^*)_{0\leq j \leq m}$ such that
$\alpha_j^* \in [a\, , \, a']$ and $a<a'$ for all $j=0,\ldots,m$,
then if $\displaystyle~~ a_N^{1+2(a'-a)} (N \, \delta_N)^{-1}
\limiteNN 0$, for all $(v_n)_n$ satisfying $\displaystyle~~ v_N
\cdot a_N^{1+2(a'-a)} (N \, \delta_N)^{-1} \limiteNN 0$,
\begin{eqnarray}\label{conv_tau}
\Prob \Big (v_N \|\underline{\tau}^*-\underline{\widehat \tau}\|_m
\geq \eta \Big ) \limiteNN 0~~~\mbox{for all $\eta>0$}.
\end{eqnarray}
\end{theorem}
Several examples of applications of this theorem will be seen in
Section \ref{applications}.
\subsection{Estimation of parameters $(\alpha_j^*)_{0\leq j \leq m}$ and $(\beta_j^*)_{0\leq j \leq m}$}
For $j=0,1,\ldots,m$, the log-log regression of $\big(S_{\widehat
k_j}^{\widehat k_{j+1}}(r_i a_N)\big)_{1\leq i\leq\ell}$ onto $(r_i
a_N)_{1\leq i\leq\ell}$ provides the estimators of $\alpha_j^*$ and
$\beta_j^*$. However, even if $ \tau_j$ converges to $\tau_j^*$,
$\widehat k_j=N\delta_N\cdot \widehat \tau_j$ does not converge to
$k_j^*$ (except if $N=o(v_N)$ which is quite impossible), and
therefore $\Prob\big([\widehat k_j,\widehat
k_{j+1}]\subset[k^*_j,k^*_{j+1}]\big)$ does not tend to $1$. So, for
$j=0,1,\ldots,m$, define $\tilde k_j$ and $\tilde k_{j}'$ such that
$$
\tilde k_j=\widehat k_j +\frac {N\delta_N}
{v_N}~~~\mbox{and}~~~\tilde k_{j}'=\widehat k_{j+1}-\frac
{N\delta_N} {v_N}~~\Longrightarrow ~~ \Prob\big([\tilde k_j,\tilde
k_{j}']\subset[k^*_j,k^*_{j+1}]\big) \limiteNN 1,
$$
from (\ref{conv_tau}) with $\eta=1/2$. Let $\Theta^*_{j}:=\left(
\begin{array}{c}
\alpha^*_j \\
\log \beta^*_j\\
\end{array} \right)$ and $\tilde\Theta_{j}:=(L_{1}^{'} \cdot L_{1})^{-1} L_{1}^{'}\cdot
Y_{\tilde k_j}^{\tilde k_{j}'}:=\left(
\begin{array}{c}
\tilde \alpha_j \\
\log \tilde \beta_j\\
\end{array} \right)$. Thus, estimators $(\tilde
\alpha_j)_{0\leq j \leq m}$ and $(\tilde \beta_j)_{0\leq j \leq m}$
satisfy
\begin{theorem}\label{convparam}
Under the same assumptions as in Theorem \ref{convprob}, for
$j=0,\ldots,m$
\begin{eqnarray}\label{conv_theta}
\sqrt {\frac {\delta_N \, N \big(\tau^*_{j+1}-\tau^*_j\big)}{a_N}}
\Big(\tilde\Theta_{j}-\Theta^*_{j}\Big)\limitedistl
\Nor\big(0,\Sigma^{(j)}(\alpha_j^*,r_1,\ldots,r_\ell)\big)
\end{eqnarray}
with $\Sigma^{(j)}(\alpha_j^*,r_1,\ldots,r_\ell):=(L_{1}^{'} \cdot
L_{1})^{-1} L_{1}^{'}\cdot\Gamma^{(j)}(\alpha_j^*,r_1,\ldots,r_\ell)
\cdot L_{1}\cdot (L_{1}^{'} \cdot L_{1})^{-1}$.
\end{theorem}
A second estimator of $\Theta^*_{j}$ can be obtained from feasible
generalized least squares (FGLS). Indeed, the asymptotic covariance
matrix $\Gamma^{(j)}(\alpha_j^*,r_1,\ldots,r_\ell)$ can be estimated
with the matrix $\tilde \Gamma^{(j)}:= \Gamma^{(j)}(\tilde
\alpha_j,r_1,\ldots,r_\ell)$ and $ \tilde \Gamma^{(j)}\limiteprob
\Gamma^{(j)}(\alpha_j^*,r_1,\ldots,r_\ell)$ since $\alpha \mapsto
\Gamma^{(j)}( \alpha,r_1,\ldots,r_\ell)$ is supposed to be a
continuous function and $\tilde \alpha_j\limiteprob\alpha^*_j$.
Since also $\alpha \mapsto \Gamma^{(j)}( \alpha,r_1,\ldots,r_\ell)$
is supposed to be a positive matrix for all $\alpha$ then
$$
\Big(\tilde
\Gamma^{(j)}\Big)^{-1}\limiteprob\Big(\Gamma^{(j)}(\alpha_j^*,r_1,\ldots,r_\ell)\Big)^{-1}
.
$$
Then, the FGLS estimator ${\overline {\Theta}}_j$ of $\Theta^*_{j}$
is defined from the minimization for all $\Theta $ of the following
criterion
$$
\parallel Y_{\tilde k_j}^{\tilde k_{j}'}-L_{a_N}\cdot \Theta
\parallel^2_{\tilde \Gamma^{(j)}}=\big(Y_{\tilde k_j}^{\tilde
k'_{j}}-L_{a_N}\cdot \Theta\big)'\cdot  \big(\tilde
\Gamma^{(j)}\big)^{-1}\cdot \big(Y_{\tilde k_j}^{\tilde
k_{j}'}-L_{a_N} \cdot \Theta\big).
$$
and therefore
$$
{\overline {\Theta}}_j=\big(L_{1}'\cdot \big(\tilde
\Gamma^{(j)}\big)^{-1}\cdot L_{1} \big)^{-1}\cdot L_{1}' \cdot
\big(\tilde \Gamma^{(j)}\big)^{-1}\cdot Y_{\tilde k_j}^{\tilde
k_{j}'}.
$$
\begin{proposition}\label{FGLS}
Under the same assumptions as in Theorem \ref{convparam}, for
$j=0,\ldots,m$
\begin{eqnarray} \label{estim2}
\sqrt {\frac {\delta_N\, N \big(\tau^*_{j+1}-\tau^*_j\big)}{a_N}}
\Big({\overline {\Theta}}_j-\Theta^*_{j}\Big)\limitedistl
\Nor\big(0,M^{(j)}(\alpha_j^*,r_1,\ldots,r_\ell)\big)
\end{eqnarray}
with $M^{(j)}(\alpha_j^*,r_1,\ldots,r_\ell):=\big(L_{1}^{'}\cdot
\big(\Gamma^{(j)}(\alpha_j^*,r_1,\ldots,r_\ell)\big)^{-1}\cdot
L_{1}\big)^{-1}\leq \Sigma^{(j)}(\alpha_j^*,r_1,\ldots,r_\ell)$
(with order's relation between positive symmetric matrix).
\end{proposition}
Therefore, the estimator ${\overline {\Theta}}_j$ converges
asymptotically faster than $\tilde \Theta_j$; ${\overline
{\alpha}}_j$ is more interesting than $\tilde \alpha_j$ for
estimating $\alpha_j^*$ when $N$ is large enough. Moreover,
confidence intervals can be easily deduced for both the estimators
of $\Theta^*_{j}$.
\subsection{Goodness-of-fit test}
For $j=0,\ldots,m$, let $T^{(j)}$ be the FGLS distance between both
the estimators of $L_{a_N} \cdot \Theta^*_{j}$, {\it i.e.} the FGLS
distance between points $\Big(\log(r_i\, a_N),\log\big(S^{\tilde
k_{j}'}_{\tilde k_j}\big)\Big)_{1\leq i\leq\ell}$ and the FGLS
regression line. The following limit theorem can be established:
\begin{theorem}\label{T}
Under the same assumptions as in Theorem \ref{convprob}, for
$j=0,\ldots,m$
\begin{eqnarray}\label{test}
T^{(j)}= \frac {\delta_N \, N \big(\tau^*_{j+1}-\tau^*_j\big)}{a_N}
\parallel Y_{\tilde k_j}^{\tilde k_{j}'}-L_{a_N}\cdot {\overline {\Theta}}_j\parallel^2_{\tilde
\Gamma^{(j)}}\limitedistl \chi^2(\ell-2).
\end{eqnarray}
\end{theorem}
{\it Mutatis mutandis}, proofs of Proposition \ref{FGLS} and Theorem
\ref{T} are the same as the proof of Proposition 5 in
\cite{bardet1}. This test can be applied to each segment $[\tilde
k_j,\tilde k_{j}'[$. However, under the assumptions, it is not
possible to prove that a test based on the sum of $T^{(j)}$ for
$j=0,\ldots,m$ converges to a $\chi^2\big ((m+1)(\ell-2)\big )$
distribution (indeed, nothing is known about the eventual
correlation of $\big(Y_{\tilde k_j}^{\tilde k_{j}'}\big)_{0\leq
j\leq m}$).
\subsection{Cases of polynomial trended processes}
Wavelet based estimators are also known to be robust to smooth
trends (see for instance \cite{Abry}). More precisely, assume now
that one considers the process $Y=\{Y_t,t\in T\}$ satisfying
$Y_t=X_t+P(t)$ for all $t \in T$ where $P$ is an unknown polynomial
function of degree $p \in \Inte$. Then,
\begin{cor}\label{trend}
Under the same assumptions as in Theorem \ref{convprob} for the
process $X$, and if the mother wavelet $\psi$ is such that $\int t^r
\, \psi(t)dt=0$ for $r=0,1,\ldots,p$, then limit theorems
(\ref{TLCS}), (\ref{lim_epsi}) and (\ref{TLClogS}) hold for $X$ and
for $Y$.
\end{cor}
Let us remark that Lemari\'e-Meyer wavelet is such that $\int t^r \,
\psi(t)dt=0$ for all $r\in \Inte$. Therefore, even if the degree $p$
is unknown, Corollary \ref{trend} can be applied. It is such the
case for locally fractional Brownian motions and applications to
heartbeat time series.
\section{Applications} \label{applications}
In this section, applications of the limit theorems to three kinds
of piecewise Gaussian processes and HR data are studied. Several
simulations for each kind of process are presented. In each case
estimators $(\widehat \tau_j)_{j}$ and $(\tilde \alpha_j)_{j}$ are
computed. To avoid an overload of results, FGLS estimators
$({\overline {\alpha}}_j)_j$ which are proved to be a little more
accurate than $(\tilde \alpha_j)_{j}$ are only presented in one case
(see Table \ref{Tablefgnbar}) because the results for $({\overline
{\alpha}}_j)_j$ are very similar to $(\tilde \alpha_j)_{j}$ ones but
are much more time consuming. For the choice of the number of scales
$\ell$, we have chosen a number proportional to the length of data
($0.15$ percent of $N$ which seems to be optimal from numerical
simulations) except in two cases (the case of goodness-of-fit test
simulations for piecewise fractional Gaussian noise and the case of
HR data, for which the length of data and the employed wavelet are
too much time consuming).
\subsection{Detection of change for Gaussian piecewise long memory
processes} In the sequel the process $X$ is supposed to be a
piecewise long range dependence time series (and therefore
$\delta_N=1$ for all $N \in \Inte$). First, some notations have to
be provided. For $Y=(Y_t)_{t\in \Inte}$ a Gaussian zero mean
stationary process, with  $r(t)=\Esp (Y_0\cdot Y_t)$ for $t\in
\Inte$, denote (when it exists) the spectral density $f$ of $Y$ by
$$
f(\lambda)=\frac 1 {2\pi} \cdot \sum_{k \in \Inter} r(k) \cdot e^{-i
k\lambda}~~\mbox{for $\lambda \in \Lambda \subset [-\pi,\pi]$}.
$$
In the sequel, the spectral density of $Y$ is supposed to satisfy
the asymptotic property,
$$
f(\lambda) \sim C \cdot \frac 1 {\lambda^{D}}~~~\mbox{when $\lambda
\to 0 $},
$$
with $C>0$ and $D\in (0,1)$. Then the process $Y$ is said to be a
long memory process and its Hurst parameter is $H=(1+D)/2$. More
precisely the following semi-parametric
framework will be considered:\\
~\\
{\bf Assumption LRD$(D)$:} $Y$ is a zero mean stationary Gaussian
process with spectral density satisfying
$$
f(\lambda)= |\lambda |^{-D} \cdot f^*(\lambda)~~\mbox{for all}~~
\lambda \in [-\pi,0[\cup ]0,\pi],
$$
with $f^*(0)>0$ and $f^*$ is such that $|f^*(\lambda)-f^*(0)|\leq
C_{2} \cdot |\lambda|^{2}~~\mbox{for all}~\lambda \in [-\pi,\pi]$
with
$C_2>0$. \\
~\\
Such assumption has been considered in numerous previous works
concerning the estimation of the long range parameter in a
semi-parametric framework (see for instance Robinson, 1995,,
Giraitis {\em et al.}, 1997, Moulines and Soulier, 2003). First and
famous examples of processes satisfying Assumption LRD$(D)$ are
fractional Gaussian noises (FGN) constituted by the increments of
the fractional Brownian motion process (FBM) and the fractionally
autoregressive integrated moving average
FARIMA$[p,d,q]$ (see more details and examples in Doukhan {\it et al.} \cite{douk}).\\
~\\
In this section, $X=(X_t)_{t\in \Inte}$ is supposed to be a Gaussian
piecewise long-range dependent process, {\it i.e.}
\begin{itemize}
\item there exists a family $(D_j^*)_{0\leq j\leq m} \in (0,1)^{m+1}$;
\item for all $j=0,\ldots,m$, for all $k \in \big
\{[N\tau_j^*],[N\tau_j^*]+1,\ldots,[N\tau_{j+1}^*]-1\big \}$,
$X_k=X^{(j)}_{k-[N\tau_{j}^*]}$ and $X^{(j)}=(X^{(j)}_t)_{t\in
\Inte}$ satisfies Assumption LRD$(D_j^*)$.
\end{itemize}
Several authors have studied the semi-parametric estimation of the
parameter $D$ using a wavelet analysis. This method has been
numerically developed by Abry {\em et al.} (1998, 2003) and Veitch
{\em et al.} (2004) and asymptotic results are provided in Bardet
{\em et al.} (2000) and recently in Moulines {\em et al.} (2007) and
Bardet {\em et al.} (2007). The following results have been
developed in this last paper. The "mother" wavelet $\psi$ is
supposed to satisfy the following assumption: first $\psi$ is
included in a Sobolev space
and secondly $\psi$ satisfies the admissibility condition. \\
~\\ {\bf  Assumption $W_1$~:} $\psi:~\Real \mapsto \Real$ with
$[0,1]$-support with $\psi(0)=\psi(1)=0$ and $ \int_0^1
\psi(t)\,dt=0$ and such that there exists a sequence
$(\psi_\ell)_{\ell \in \Inter}$ satisfying $\psi(\lambda)=\sum_{\ell
\in \Inter}\psi_\ell e^{2\pi i\ell \lambda} \in\LL^2([0,1])$ and
$\sum_{\ell\in\Inter}(1+|\ell|)^{5/2} |\psi_\ell|<\infty$.\\
~\\
For ease of writing, $\psi$ is supposed to be supported in $[0,1]$.
By an easy extension the following propositions are still true for
any compactly supported wavelets. For instance, $\psi$ can be a
dilated Daubechies "mother" wavelet of order $d$ with $d \geq 6$ to
ensure the smoothness of the function $\psi$. However, the following
proposition could also be extended for "essentially" compactly
supported "mother" wavelet like Lemari\'e-Meyer wavelet. Remark that
it is not necessary to choose $\psi$ being a "mother" wavelet
associated to a multi-resolution analysis of $\LL^2(\Real)$ like in
the recent paper of Moulines {\em et al.} (2007). The whole theory
can be developed without resorting to this assumption. The choice of
$\psi$ is then very large. Then, in Bardet {\it et al.} (2007), it
was established:
\begin{propo}\label{tlclog}
Let $X$ be a Gaussian piecewise long-range dependent process defined
as above and $(a_n)_{n \in \Inte}$ be such that $N/a_N \limiteNN
\infty$ and $a_N \cdot N^{-1/5}\limiteNN \infty$. Under Assumption
$W_1$, limit theorems (\ref{TLCS}), (\ref{lim_epsi}) and
(\ref{TLClogS}) hold with $\alpha_j^*=D_j^*$ and $\beta_j^*=\log
\Big (f_j^*(0)\, \int_{-\infty}^\infty |\widehat \psi (u)|^2 \cdot
|u|^{-D}du\Big )$ for all $j=0,1,\ldots,m$ and with
$d_{pq}=GCD(r_p\, , \, r_q)$ for all $(p,q)\in \{1,\ldots,\ell\}$,
\begin{eqnarray*}\label{cov}
\gamma^{(j)}_{pq}=\frac {2(r_pr_q)^{2-D_j^*}} {d_{pq}}
\sum_{m=-\infty}^{\infty}\left ( \frac {\int _0 ^\infty \widehat
\psi (ur_p)\overline{\widehat \psi} (ur_q)\, u^{-D_j^*}\,\cos (u
\,d_{pq} m)\, du}{\int_{0}^\infty |\widehat \psi (u)|^2 \cdot
|u|^{-D_j^*}du} \right) ^2.
\end{eqnarray*}
\end{propo}
As a consequence, the results of Section \ref{results} can be
applied to Gaussian piecewise long-range dependent processes:
\begin{cor}\label{cor_LRD}
Under assumptions of Proposition \ref{tlclog} and Assumption C, for
all $0<\kappa< 2/15$, if $a_N=N^{\kappa+1/5}$ and
$v_N=N^{2/5-3\kappa}$ then (\ref{conv_tau}), (\ref{conv_theta}),
 (\ref{estim2}) and (\ref{test}) hold.
\end{cor}
Thus, the rate of convergence of $\underline{\widehat \tau}$ to
$\underline{\tau}^*$ (in probability) is $N^{2/5-3\kappa}$ for
$0<\kappa$ as small as one wants. Estimators $\tilde D_j$ and
${\overline D_j}$ converge to the parameters $D_j^*$ following a
central limit theorem with a rate of convergence
$N^{2/5-\kappa/2}$ for $0<\kappa$ as small as one wants.\\
~\\
{\it Results of simulations:} The following Table \ref{Tablefgn}
represents the change point and parameter estimations in the case of
a piecewise FGN with one abrupt change point. We observe the good
consistence property of the estimators. Kolmogorov-Smirnov tests
applied to the sample of estimated parameters lead to the following
results: \begin{enumerate}
\item the estimator $\widehat \tau_1$ can not be modeled with a Gaussian distribution;
\item the estimator $\widehat H_j$ seems to follow a Gaussian
distribution.
\end{enumerate}
\begin{table}[!h]
\renewcommand{\arraystretch}{1.2}
\begin{center}
\begin{tabular} {|c|}
\hline
$N=20000$, $\tau_1=0.75$,
$D_0=0.2$ and $D_1=0.8$\\
\hline \hline
\begin{tabular} {ccc|ccc|ccc} $\widehat \tau_1$ & $\widehat \sigma_{\tau_1}$ &
$\sqrt{MSE}$ & $\tilde D_0$ & $\widehat \sigma_{D_0}$ & $\sqrt{MSE}$
& $\tilde D_1$ & $\widehat \sigma_{D_1}$ &
$\sqrt{MSE}$\\
\hline
0.7605 & 0.0437  & 0.0450 & 0.2131 & 0.0513 & 0.0529 & 0.7884 & 0.0866 & 0.0874\\
\end{tabular}\\
\hline
\end{tabular}
\end{center}
\caption{Estimation of $\tau_1$, $D_0$ and $D_1$ in the case of a
piecewise FGN ($H_0=0.6$ and $H_1=0.9$) with one change point when
$N=20000$ and $\ell=30$ (50 realizations)}\label{Tablefgn}
\end{table}
The distribution of the test statistics $T^{(0)}$ and $T^{(1)}$ (in
this case $\ell=20$ and $N=20000$ and 50 realizations) are compared
with a Chi-squared-distribution with eighteen degrees of freedom.
The goodness-of-fit Kolmogorov-Smirnov test for $T^{(j)}$ to the
$\chi^2(18)$-distribution is accepted (with $0.3459$ for the sample
of $T^{(0)}$ and $p=0.2461$ for $T^{(1)}$). In this case and for the
same parameters as in Table \ref{Tablefgn}, the estimator
${\overline {D}}_j$ seems to be a little more accurate than $\tilde
D_j$ (see Table \ref{Tablefgnbar}).
\begin{table}[!h]
\renewcommand{\arraystretch}{1.2}
\begin{center}
\begin{tabular} {|ccc|ccc|ccc|}
\hline $\widehat \tau_1$ & $\widehat \sigma_{\tau_1}$ & $\sqrt{MSE}$
& ${\overline {D}}_0$ & $\widehat \sigma_{D_0}$ & $\sqrt{MSE}$ &
${\overline {D}}_1$ & $\widehat \sigma_{D_1}$ &
$\sqrt{MSE}$\\
\hline 0.7652 & 0.0492  & 0.0515 & 0.1815 & 0.0452 & 0.0488 & 0.8019
& 0.0721 & 0.0722\\
\hline
\end{tabular}
\end{center}
\caption{Estimation of $D_0$ and $D_1$ in the case of a piecewise
FGN ($D_0=0.2$ and $D_1=0.8$) with one change point when $N=20000$
and $\ell=20$ (50 realizations)}\label{Tablefgnbar}
\end{table}\\
~\\
Simulations are also applied to a piecewise simulated
FARIMA(0,$d_j$,0) processes and results are similar (see Table
\ref{Tablefarima}). The following Figure \ref{FARIMA} represents the
change point instant and its estimation for such a process with one
abrupt change point.
\begin{center}
\begin{figure}[htbp]
\begin{center}
\includegraphics[width=10 cm,height=4 cm]{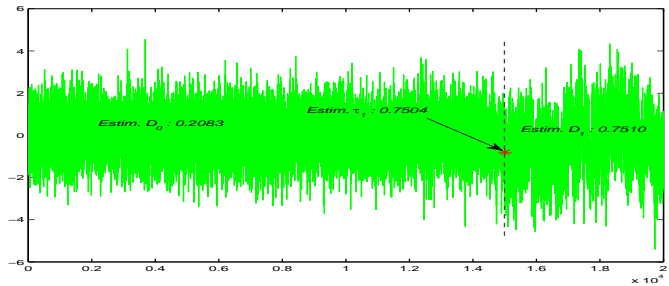}
\end{center}
\caption{Detection of the change point in piecewise
FARIMA(0,$d_j$,0) (for the first segment $d_0=0.1$ ($D_0=0.2$) for
the second $d_1=0.4$ ($D_1=0.8$))} \label{FARIMA}
\end{figure}
\end{center}
\begin{table}[!h]
\renewcommand{\arraystretch}{1.2}
\begin{center}
\begin{tabular} {|c|}
\hline $N=20000$, $\tau_1=0.75$,
$D_0=0.2$ and $D_1=0.8$ \\
\hline \hline
\begin{tabular} {ccc|ccc|ccc}$\widehat \tau_1$ & $\widehat \sigma_{\tau_1}$ &
$\sqrt{MSE}$ & $\tilde D_0$ & $\widehat \sigma_{D_0}$ & $\sqrt{MSE}$
& $\tilde D_1$ & $\widehat \sigma_{D_1}$ &
$\sqrt{MSE}$\\
\hline 0.7540 & 0.0215  & 0.0218  & 0.1902 & 0.0489  & 0.0499 & 0.7926 & 0.0761 & 0.0764\\
\end{tabular}\\
\hline
\end{tabular}
\end{center}
\caption{Estimation of $\tau_1$, $D_0$ and $D_1$ in the case of
piecewise FARIMA(0,$d_j$,0) ($d_0=0.1$ and $d_1=0.4$) with one
change point when $N=20000$ and $\ell=30$ ($50$
realizations)}\label{Tablefarima}
\end{table}
\subsection{Detection of abrupt change for piecewise Gaussian
self-similar processes} Let us recall that $B^H=(B^H_t)_{t \in
\Real}$ is a fractional Brownian motion (FBM) with two parameters $H
\in (0,1)$  and $\sigma^2>0$ when $B_H$ is a Gaussian process having
stationary increments and such as
$$
\Var (B^H_t))= \sigma ^{2}|t|^{2H} \,\,\,\, \forall t\in \Real.
$$
It can be proved that $B_H$ is the only Gaussian self-similar
process having stationary increments and its self-similar parameter
is $H$ (a process $Y=(Y_t)_{t \in E}$ is said to be a
$H_s$-self-similar process if for all $c>0$ and for all
$(t_1,\ldots,t_k)\in E^k$ where $k\in \Inte^*$, the vector
$\big(Y_{ct_1},\ldots,Y_{ct_k}\big )$ has the same distribution
than the vector $c^{H_s}\big(Y_{t_1},\ldots,Y_{t_k}\big )$).\\
~\\
Now, $X$ will be called a piecewise fractional Brownian motion if:
\begin{itemize}
\item there exist two families of parameters $(H^*_j)_{0\leq j\leq m}\in (0,1)^{m+1}$
and $(\sigma^{*2}_j)_{0\leq j\leq m}\in (0,\infty)^{m+1}$;
\item for all $j=0,\ldots,m$, for all $t \in \big
[[N\tau_j^*],[N\tau_j^*]+1,\ldots,[N\tau_{j+1}^*]-1\big ]$,
$X_t=X^{(j)}_{t-[N\tau_{j}^*]}$ and $X^{(j)}=(X^{(j)}_t)_{t\in
\Real}$ is a FBM with parameters $H^*_j$ and $\sigma^{*2}_j$.
\end{itemize}
The wavelet analysis of FBM has been first studied by Flandrin
(1992) and developed by Abry (1998) and Bardet (2002). Following
this last paper, the mother wavelet $\psi$ is supposed to
satisfy:\\
~\\
{\bf Assumption $W_2$}: $\psi:\Real \to \Real$ is a piecewise
continuous and left (or right)-differentiable in $[0,1]$, such that
$|\psi'(t^-)|$ is Riemann integrable in $[0,1]$ with $\psi'(t^-)$
the left-derivative of $\psi$ in $t$, with support included in
$[0,1]$ and $ \int_\Real t^p\psi(t)\,dt=\int_0^1
t^p\psi(t)\,dt=0~~\mbox{for}~~p=0,1.$\\
~\\
As in Assumption $W_1$, $\psi$ is supposed to be supported in
$[0,1]$ but the following propositions are still true for any
compactly supported wavelets. Assumption $W_2$ is clearly weaker
than Assumption $W_1$ concerning the regularity of the mother
wavelet. For instance, $\psi$ can be a Daubechies wavelet of order
$d$ with $d \geq 3$ (the Haar wavelet, {\it i.e.} $d=2$, does not
satisfy $\int_0^1 t\,\psi(t)\,dt=0$). Another choice could be
infinite support wavelets with compact effective support (it is such
the case with Meyer or Mexican Hat wavelets) but the proof of the
following property has to be completed.
\begin{propo} \label {TLCprop}
Assume that $X$ is a piecewise FBM as it is defined above and let
$(X_1,X_2,\ldots,X_N)$ be a sample of a path of $X$ (therefore
$\delta_N=1$). Under Assumption $W_2$, if $(a_n)_{n \in \Inte}$ is
such that $N/a_N \limiteNN \infty$ and $a_N \cdot N^{-1/3}\limiteNN
\infty$, then limit theorems (\ref{TLCS}), (\ref{lim_epsi}) and
(\ref{TLClogS}) hold with $\alpha_j^*=2H_j^*+1$ and $\beta_j^*=\log
\big (-\frac {\sigma_j^{*2}} 2 \int_0^1 \int_0^1 \psi(t) \psi(t') |
t-t' | ^{2H_j^*}dt\, dt'\big)$ for all $j=0,1,\ldots,m$ and with
$d_{pq}=GCD(r_p\, , \, r_q)$ for all $(p,q)\in \{1,\ldots,\ell\}$,
\begin{eqnarray*} \label{F2}
\gamma^{(j)}_{pq}= \frac {2d_{pq}} {r_p^{2H_j^*+1/2}r_q^{2H_j^*+1/2}
} \sum_{k=-\infty}^{\infty}\left ( \frac {\int _0^1 \int _0^1
\psi(t) \psi(t') \left | k\,d_{pq}+r_pt-r_qt' \right | ^{2H_j^*}
dt\, dt'}{\int_0^1 \int_0^1 \psi(t) \psi(t') | t-t' | ^{2H_j^*}dt\,
dt'}\right )^2 .
\end{eqnarray*}
\end{propo}
Then, Theorem \ref{convprob} can be applied to piecewise FBM but
$2(a'-a)+1=2(\sup_j\alpha_j^*-\inf_j\alpha_j^*)+1$ has to be smaller
than $3$ since $a_N \cdot N^{-1/3}\limiteNN \infty$. Thus,
\begin{cor}\label{cor_FBM}
Let $A:=\big |\sup_jH_j^*-\inf_jH_j^*\big |$. If $A<1/2$, under
assumptions of Proposition \ref{TLCprop} and Assumption C, for all
$0<\kappa<\frac 1 {1+4A}-\frac 1 3$, if $a_N=N^{1/3+\kappa}$ and
$v_N=N^{2/3(1-2A)-\kappa(2+4A)}$ then (\ref{conv_tau}),
(\ref{conv_theta}), (\ref{estim2}) and (\ref{test}) hold.
\end{cor}
Thus, the rate of convergence of $\underline{\widehat \tau}$ to
$\underline{\tau}^*$ (in probability) can be $N^{2/3(1-2A)-\kappa'}$
for $0<\kappa'$ as small as one wants when
$a_N=N^{1/3+\kappa'/(2+4A)}$. \\
~\\
{\bf Remark:} This result of Corollary \ref{cor_FBM} is quite
surprising: the smaller $A$, {\it i.e.} the smaller the differences
between the parameters $H_j$, the faster the convergence rates of
estimators $\widehat \tau_j$ to $\tau_j^*$. And if the difference
between two successive parameters $H_j$ is too large, the estimators
${\widehat\tau_j}$ do not seem to converge. Following simulations in
Table \ref{Tableecart} will exhibit this paroxysm. This induces a
limitation of the estimators' using especially for applying them to
real data (for
which a priori knowledge is not available about the values of $H_j^*$). \\
~\\
Estimators $\tilde H_j$ and ${\overline H_j}$ converge to the
parameters $H_j^*$ following a central limit theorem with a rate of
convergence $N^{1/3-\kappa/2}$ for $0<\kappa$ as small as one
wants. \\
~\\
{\it Results of simulations:} The following Table \ref{Tablefbm1}
represent the change point and parameter estimations in the case of
piecewise FBM with one abrupt change point. Estimators of the change
points and parameters seem to converge since their mean square
errors clearly decrease when we double the number of observations.
\begin{table}[!h]
\renewcommand{\arraystretch}{1.2}
\begin{center}
\begin{tabular} {|c||c||c|}
\hline & $N=5000$ & $N=10000$\\
\hline
\begin{tabular} {c}
$\tau_1$ \\ 0.4  \\ \hline $H_0$ \\ 0.4 \\ \hline $H_1$ \\0.8
\end{tabular} &
\begin{tabular} {c|c|c}
$\widehat \tau_1$ & $\widehat \sigma_{\tau_1}$ &
$\sqrt{MSE}$ \\
0.4467  & 0.0701 & 0.0843 \\
\hline
$\tilde H_0$ & $\widehat \sigma_{H_0}$ & $\sqrt{MSE}$ \\
0.3147  & 0.0404 & 0.0943  \\
\hline $\tilde H_1$ & $\widehat \sigma_{H_1}$ & $\sqrt{MSE}$ \\
0.7637  & 0.0534 & 0.0645
\end{tabular}
&
\begin{tabular} {c|c|c}
$\widehat \tau_1$ & $\widehat \sigma_{\tau_1}$ &
$\sqrt{MSE}$ \\
0.4368 & 0.0319 & 0.0487 \\
\hline $\tilde H_0$ & $\widehat \sigma_{H_0}$ & $\sqrt{MSE}$ \\
0.3761 & 0.0452 & 0.0511 \\
\hline $\tilde H_1$ & $\widehat \sigma_{H_1}$ & $\sqrt{MSE}$ \\
0.7928 & 0.0329 & 0.0337
\end{tabular} \\
\hline
\end{tabular}
\end{center}
\caption{Estimation of $\tau_1$, $H_0$ and $H_1$ in the case of
piecewise FBM with one change point when $N=5000$ ($100$
realizations) and $N=10000$ ($50$ realizations)} \label{Tablefbm1}
\end{table}
~\\
For testing if the estimated parameters follow a Gaussian
distribution, Kolmogorov-Smirnov goodness-of-fit tests (in the case
with $N=10000$ and $50$ replications) are applied:
\begin{enumerate}
\item this test  for $\tilde H_0$ is accepted as well as for
$\tilde H_1$ and the following Figure \ref{DENS} represents the
relating distribution.
\item this is not such the case for the
change point estimator $\widehat \tau_1$ for which the hypothesis of
a possible fit with a Gaussian distribution is rejected
($KS_{test}=0.2409$) as showed in the Figure \ref{DENS} below which
represents the empirical distribution function with the
correspondant Gaussian cumulative distribution function.
\end{enumerate}
\begin{center}
\begin{figure}[h]
\begin{center}
\includegraphics[width=7 cm,height=5 cm]{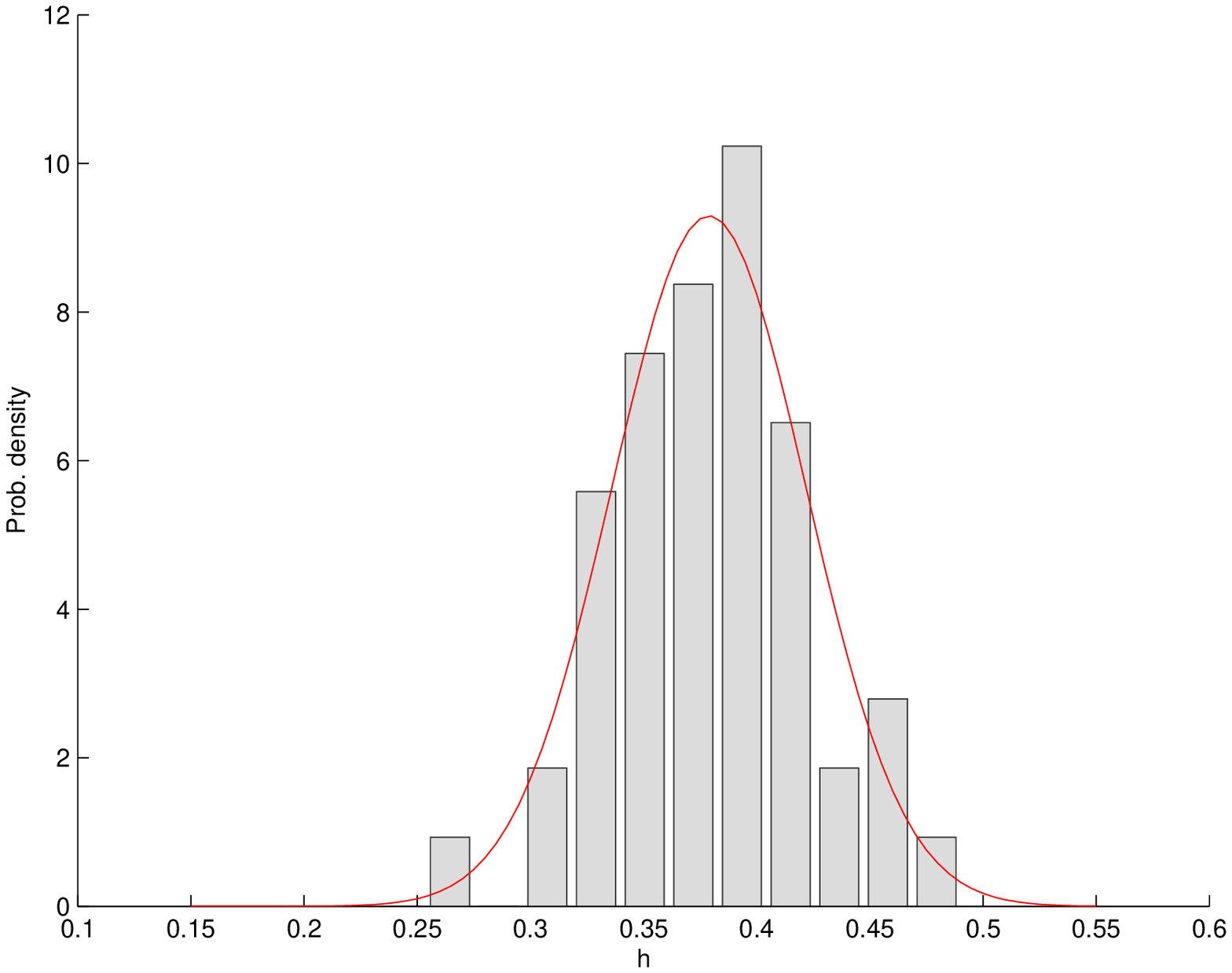}
\vspace{1cm}\includegraphics[width=7 cm,height=5 cm]{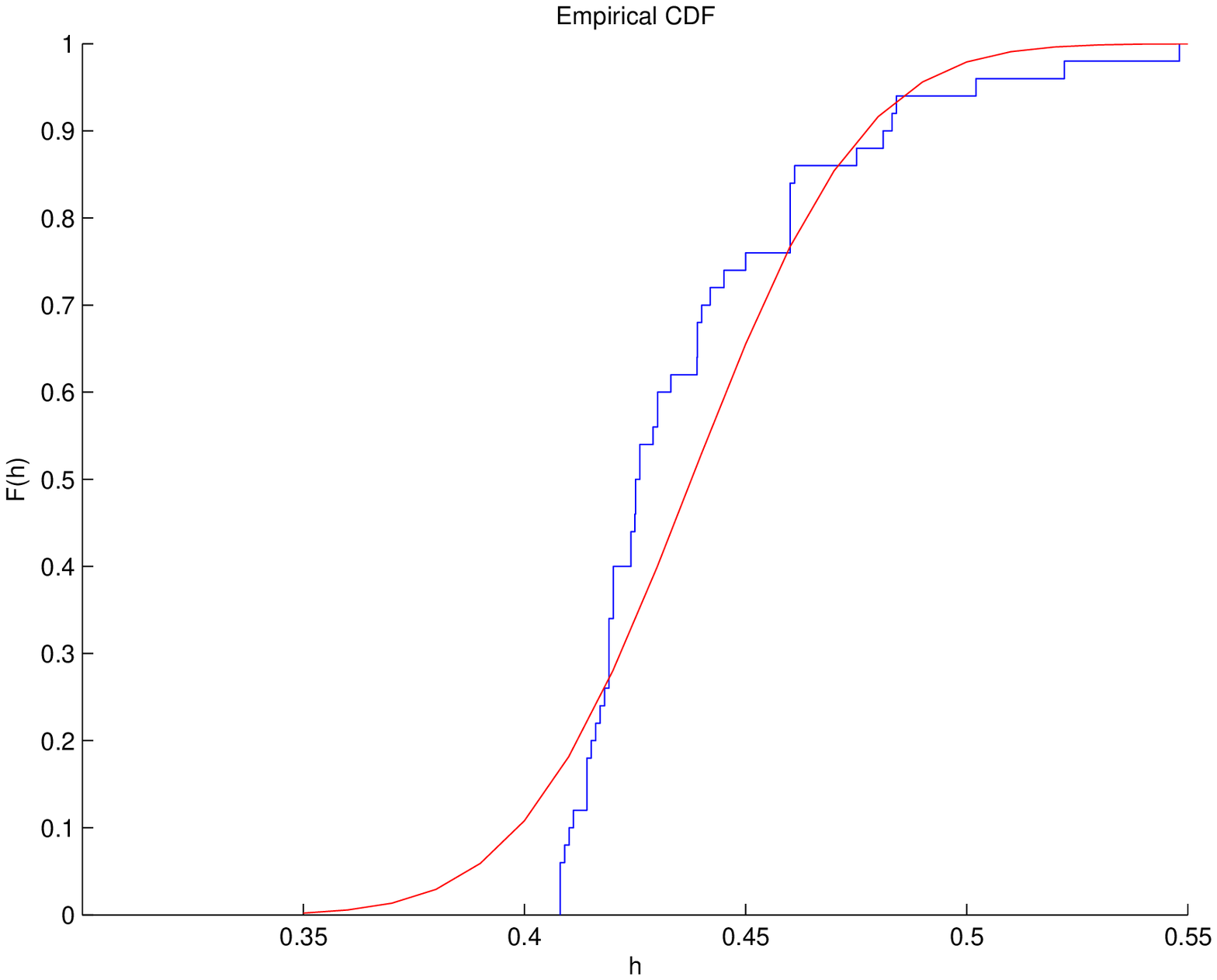}
\end{center}
\caption{Left: Modeling of sample estimations of $\tilde H_0$ with
normal distribution; Right: Comparison of the generated empirical
cumulative distribution for $\widehat \tau_1$ (when N=10000) and the
theoretical normal distribution.} \label{DENS}
\end{figure}
\end{center}
From the following example in Table \ref{Tableecart}, we remark that
the estimated parameters seem to be non convergent when the
difference between the parameters $H_j$ is too large.
\begin{table}[!h]
\renewcommand{\arraystretch}{1.2}
\begin{center}
\begin{tabular} {|c|}
\hline $N=5000$, $\tau_1=0.6$,
$H_0=0.1$ and $H_1=0.9$ \\
\hline \hline
\begin{tabular} {ccc|ccc|ccc}
$\widehat \tau_1$ & $\widehat \sigma_{\tau_1}$ & $\sqrt{MSE}$ &
$\tilde H_0$ & $\widehat \sigma_{H_0}$ & $\sqrt{MSE}$ & $\tilde H_1$
& $\widehat \sigma_{H_1}$ &
$\sqrt{MSE}$\\
\hline 0.5950 & 0.1866  & 0.1866 & -0.1335 & 0.0226  & 0.2346 &
0.6268 & 0.4061 & 0.4894
\end{tabular}\\
\hline
\end{tabular}
\end{center}
\caption{Estimation of $\tau_1$, $H_0$ and $H_1$ (when
$H_1-H_0=0.8>1/2$) in the case of piecewise FBM with one change
point when $N=5000$ (50 realizations)}\label{Tableecart}
\end{table}
~\\
Simulations for goodness-of-fit tests $T^{(j)}$ provide the
following results: when $N=5000$, the drawn distributions of the
computed test statistics (see Figure \ref{plotkhi2}) exhibit a
Khi-square distributed values ($\chi^2(5)$ since $\ell=7$) and $95
\%$ of the $100$ of the values of $T^{(0)}$ and $T^{(1)}$ do not
exceed $\chi^2_{95\%}(5)=11.0705$. These results are also validated
with Kolmogorov-Smirnov tests.
\begin{figure}[h!]
\begin{minipage}[b]{.5\linewidth}
\begin{center}
\includegraphics[width=5 cm,height=4 cm]{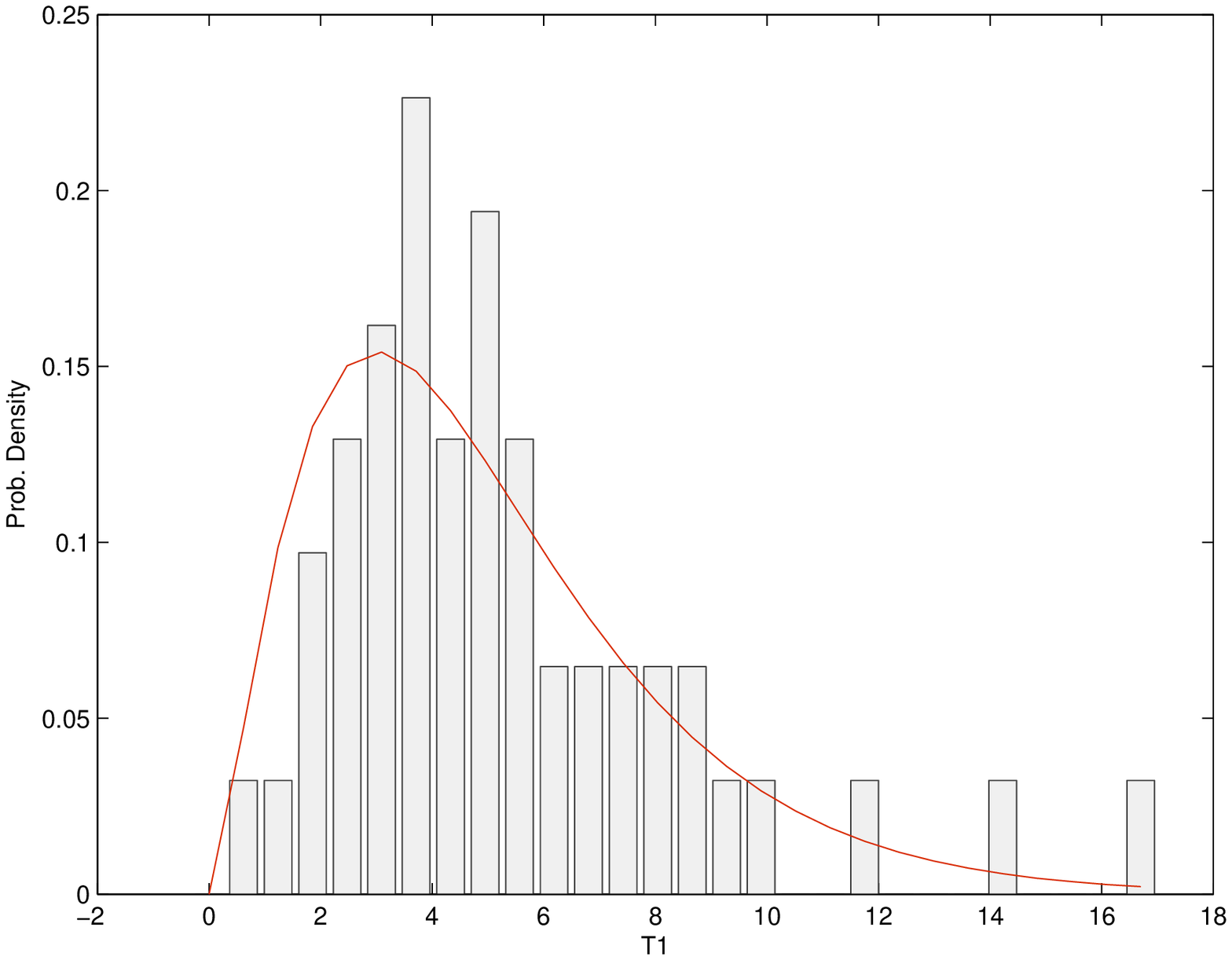}
\end{center}
\end{minipage}
\begin{minipage}[b]{.4\linewidth}
\begin{center}
\includegraphics[width=5 cm,height=4 cm]{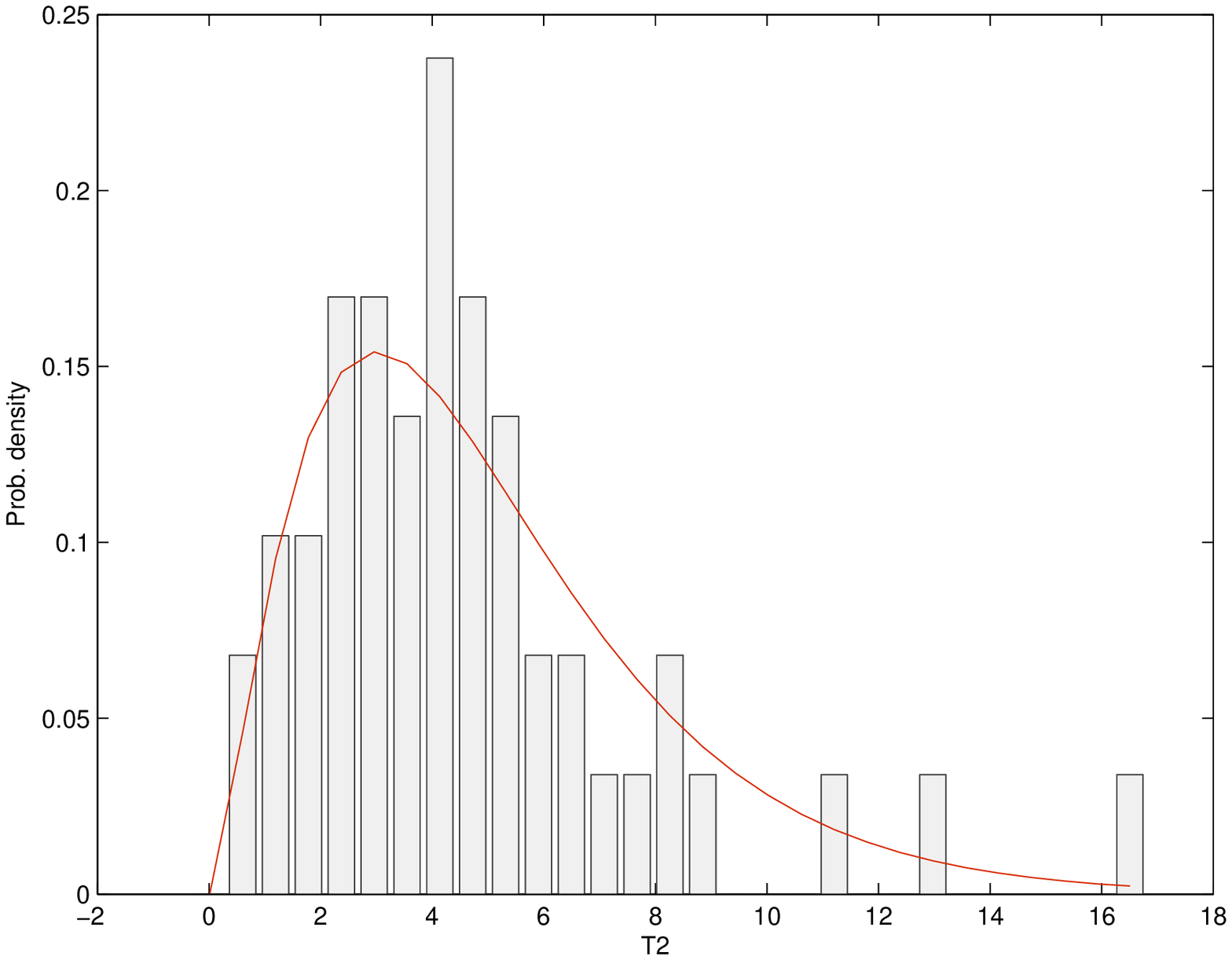}
\end{center}
\end{minipage}
\caption{Testing for $\chi^2(5)$ distribution in the first detected
zone (left) and the second detected zone (right) ($50$ realizations
when $N=5000$)} \label{plotkhi2}
\end{figure}
\newline The results below in Table \ref{Tablefbm2} are obtained with
piecewise fractional Brownian motion when two change points are
considered. As previously, both the $KS_{test}$ tests for deciding
whether or not samples of both estimated change points is consistent
with Gaussian distributions are rejected. However, such $KS_{test}$
tests are accepted for $\tilde H_j$ samples. A graphical
representation of the change point detection method applied to a
piecewise FBM is given in Figure \ref{Ruptfbm}.\\
\begin{table}[!h]
\renewcommand{\arraystretch}{1.2}
\begin{center}
\begin{tabular} {|c||c||c|}
\hline & $N=5000$ & $N=10000$\\
\hline
\begin{tabular} {c}
$\tau_1$ \\ 0.3  \\ \hline $\tau_2$ \\ 0.78  \\ \hline $H_0$ \\ 0.6 \\
\hline $H_1$ \\0.8 \\ \hline $H_2$ \\0.5
\end{tabular} &
\begin{tabular} {c|c|c}
$\widehat \tau_1$ & $\widehat \sigma_{\tau_1}$ &
$\sqrt{MSE}$ \\
0.3465  & 0.1212 & 0.1298 \\
\hline $\widehat \tau_2$ & $\widehat \sigma_{\tau_2}$ &
$\sqrt{MSE}$ \\
0.7942  & 0.1322 & 0.1330 \\
\hline
$\tilde H_0$ & $\widehat \sigma_{H_0}$ & $\sqrt{MSE}$ \\
0.5578  & 0.0595 & 0.0730  \\
\hline
$\tilde H_1$ & $\widehat \sigma_{H_1}$ & $\sqrt{MSE}$ \\
0.7272  & 0.0837 & 0.1110 \\
\hline
$\tilde H_2$ & $\widehat \sigma_{H_2}$ & $\sqrt{MSE}$ \\
0.4395  & 0.0643 & 0.0883
\end{tabular}
&
\begin{tabular} {c|c|c}
$\widehat \tau_1$ & $\widehat \sigma_{\tau_1}$ &
$\sqrt{MSE}$ \\
0.3086 & 0.0893 & 0.0897 \\
\hline $\widehat \tau_2$ & $\widehat \sigma_{\tau_2}$ &
$\sqrt{MSE}$ \\
0.7669 & 0.0675 & 0.0687 \\
\hline
$\tilde H_0$ & $\widehat \sigma_{H_0}$ & $\sqrt{MSE}$ \\
0.5597  & 0.0449 & 0.0604 \\
\hline
$\tilde H_1$ & $\widehat \sigma_{H_1}$ & $\sqrt{MSE}$ \\
0.7633   & 0.0813 & 0.0892 \\
\hline
$\tilde H_2$ & $\widehat \sigma_{H_2}$ & $\sqrt{MSE}$ \\
0.4993   & 0.0780 & 0.0780
\end{tabular} \\
\hline
\end{tabular}
\end{center}
\caption{Estimation of $\tau_1$, $\tau_2$, $H_0$, $H_1$ and $H_2$ in
the case of piecewise FBM with two change points when $N=5000$ and
$N=10000$ (50 realizations)} \label{Tablefbm2}
\end{table}
\begin{figure}[h!]
\begin{minipage}[c]{.6\linewidth}
\begin{center}
\includegraphics[width=8 cm,height=5 cm]{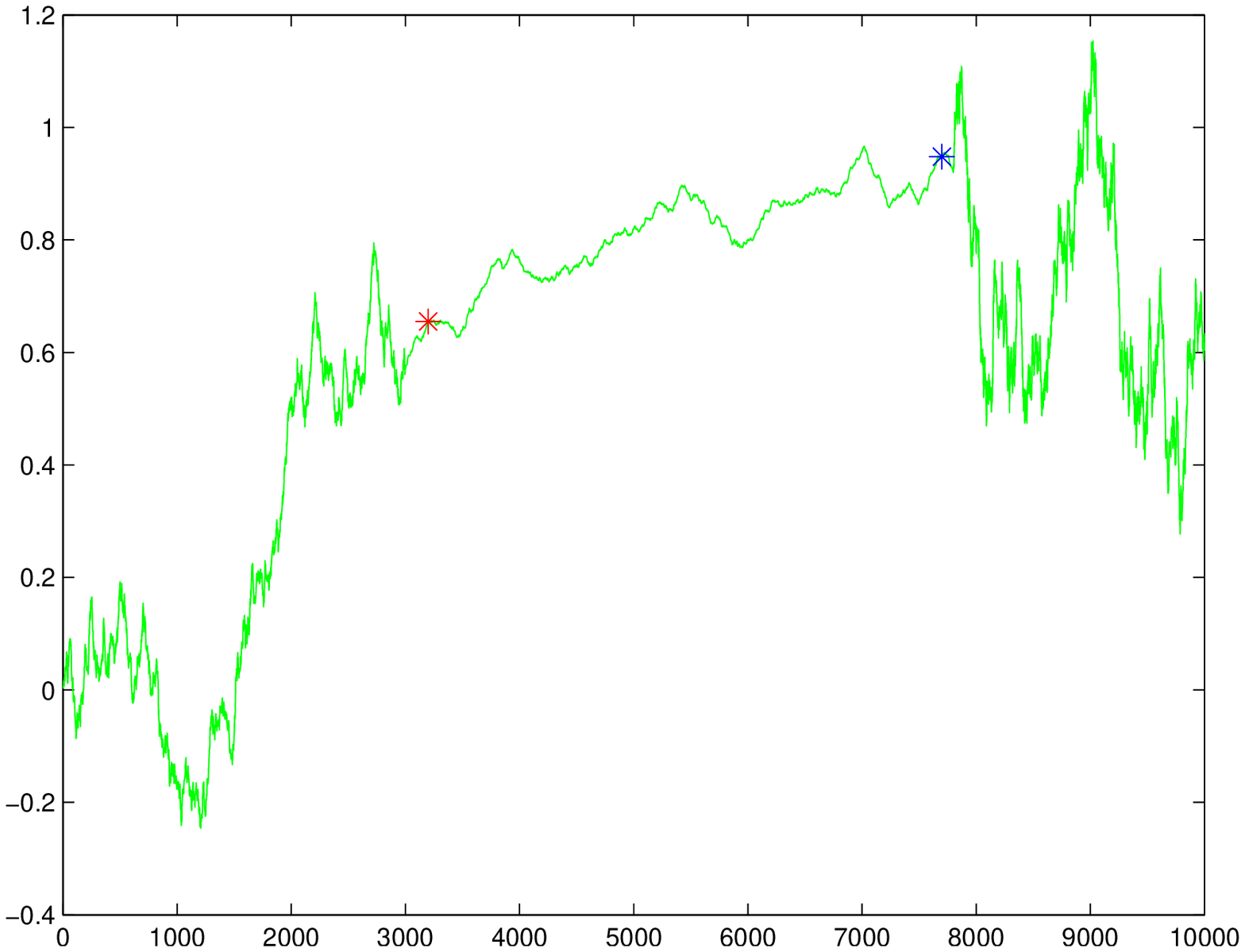}
\end{center}
\end{minipage}
\begin{minipage}[c]{.3\linewidth}
\begin{center}
\includegraphics[width=5 cm,height=4 cm]{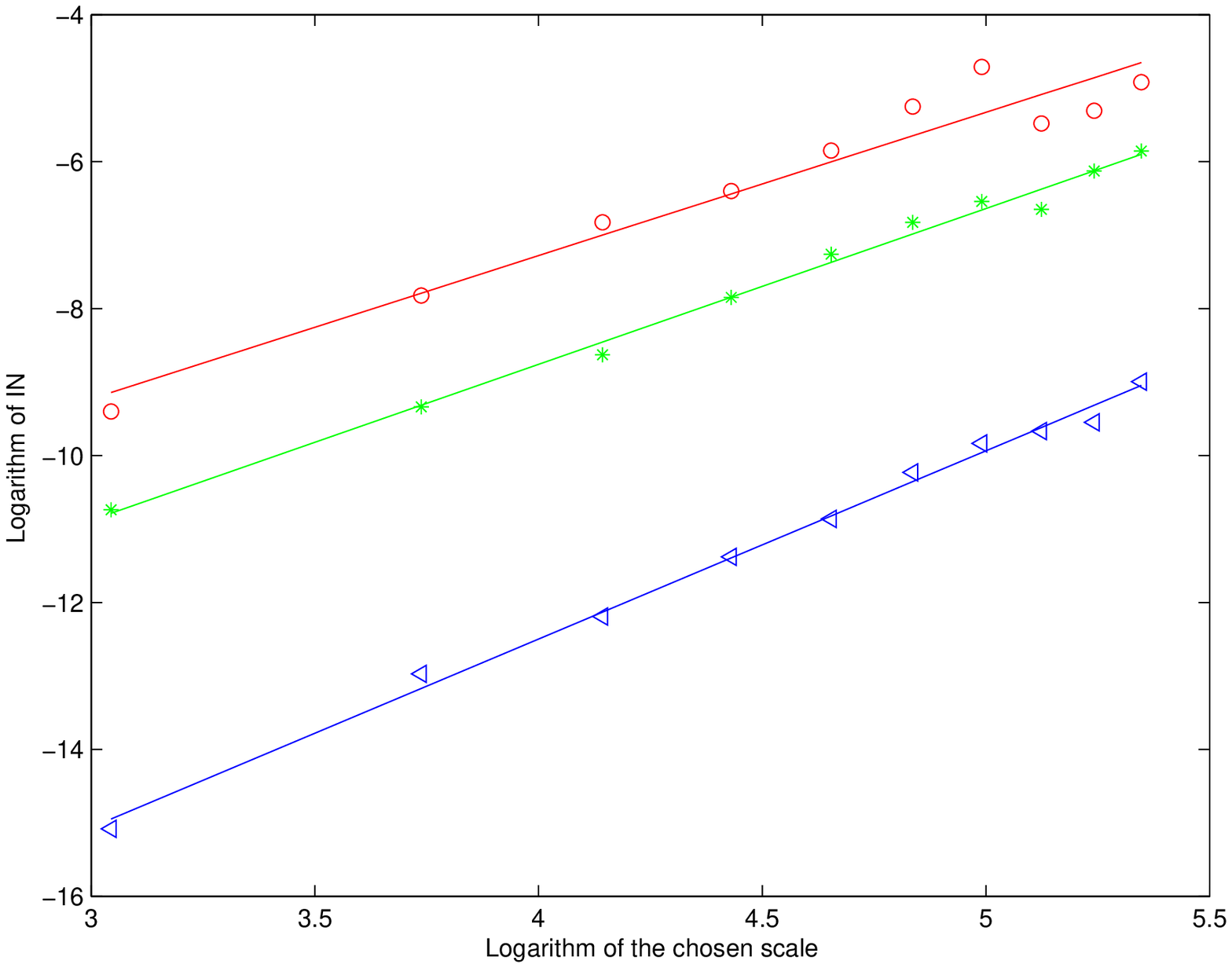}
\end{center}
\end{minipage}
\caption{(left)Detection of the change point in piecewice FBM($H_j$)
($\tau_1=0.3$, $\tau_2=0.78$, $H_0=0.6$, $H_1=0.8$ and $H_2=0.5$).
The change points estimators are $\widehat\tau_1=0.32$ and
$\widehat\tau_2=0.77$. (right) Representation of log-log regression
of the variance of wavelet coefficients on the chosen scales for the
three segments ($\tilde H_0=0.5608$ (*), $\tilde H_1=0.7814$
($\lhd$) and $\tilde H_2=0.4751$ (o)) } \label{Ruptfbm}
\end{figure}
~\\
The distribution of the test statistics $T^{(0)}$, $T^{(1)}$ and
$T^{(2)}$ (in this case $\ell=10$, $N=10000$ and 50 realizations)
are compared with a Chi-squared-distribution with eight degrees of
freedom. The goodness-of-fit Kolmogorov-Smirnov test for $T^{(j)}$
to the $\chi^2(8)$-distribution is accepted (with $p=0.4073$ for the
sample of $T^{(0)}$, $p=0.2823$ for $T^{(1)}$ and $p=0.0619$ for
$T^{(2)}$).
\subsection{Detection of abrupt change for piecewise locally fractional Gaussian
processes} In this section, a continuous-time process $X$ is
supposed to model data. Therefore assume that
$(X_{\delta_N},X_{2\delta_N},\ldots,X_{N\, \delta_N})$ is known,
with $\delta_N \limiteNN 0$ and $N \, \delta_N \limiteNN \infty$. A
piecewise locally fractional Gaussian process $X=(X_t)_{t \in
\Real_+}$ is defined by
\begin{eqnarray}\label{def_Xss}
X_t:=\int_{\Real}\frac{e^{it\xi}-1}{\rho_j(\xi)}\widehat{W}(d\xi)~~\mbox{for}~~t\in
[\tau_j^*N \, \delta_N\, , \, \tau_{j+1}^*N \, \delta_N)
\end{eqnarray}
where the functions $\rho_j: \Real \to [0,\infty)$ are even Borelian
functions such that for all $j=0,1,\ldots,m$,:
\begin{itemize}
\item $\displaystyle \rho_j(\xi)= \frac 1 {\sigma^*_j} \,
|\xi|^{H^*_j+{1/2}}~~\mbox{for}~~ |\xi|\in[f_{min}\, , \,
f_{max}]~~\mbox{with}~~H^*_j \in \Real,~\sigma^*_j>0$;
\item $\displaystyle \int_{\Real} \left(1\wedge  |\xi|^2\right) \frac 1 {\rho_j^2(\xi)}
\, d \xi < \infty$
\end{itemize}
and $W(dx)$ is a Brownian measure and $\widehat{W}(d\xi)$ its
Fourier transform in the distribution meaning. Remark that
parameters $H_j^*$, called local-fractality parameters, can be
supposed to be included in $\Real$ instead the usual interval
$(0,1)$. Here $0<f_{min}<f_{max}$ are supposed to be known
parameters. Roughly speaking, a locally fractional Gaussian process
is nearly a self-similar Gaussian process for scales (or
frequencies) included in a band of scales
(frequencies). \\
~\\
For locally fractional Gaussian process already studied in Bardet
and Bertrand (2007) and Kammoun {\it et al.} (2007), the mother
wavelet is supposed to satisfy \\
~\\
{\bf Assumption $W_3$}: $\psi:~\Real \mapsto \Real$ is a ${\cal
C}^\infty(\Real)$ function such that for all $m\in \Inte$,
$\int_{\Real} \left |t^m\psi(t)\right |dt <\infty$ and the Fourier
transform $\widehat{\psi}$ of $\psi$ is an even function compactly
supported
on $[-\mu,-\lambda]\cup[\lambda,\mu]$ with $0 <\lambda<\mu$.\\
~\\
These conditions are sufficiently mild and are satisfied in
particular by the Lemari\'e-Meyer "mother" wavelet. The
admissibility property, i.e. $\int_{\Real} \psi(t)dt = 0$, is a
consequence of the second condition and more generally, for all $m
\in \Inte$, $\int_{\Real}t^m \psi(t)dt = 0.$\\
~\\
Since the function $\psi$ is not a compactly supported mother
wavelet, wavelet coefficients $d_X(a,b)$ can not be well
approximated by $e_X(a,b)$ when the shift $b$ is close to $0$ or
$N\, \delta_N$. Then, a restriction $\tilde S^{k'}_{k}(a_N)$ of
sample wavelet coefficient's variance $S^{k'}_{k}(a_N)$ has to be
defined:
\begin{eqnarray*}\label{Sktilde}
\tilde S^{k'}_{k}(a_N):= \frac {a_N} {(1-2w)k'-k} \, \sum
_{p=[(k+w(k'-k))/a_N]+w}^{[(k'-w(k'-k))/a_N]-1}e_{X}^2(a_N,a_N\,p)~~\mbox{with
$0<w<1/2$}.
\end{eqnarray*}
\begin{propo} \label{Theo:TCL:discret}
Assume that $X$ is a piecewise locally fractional Gaussian process
as it is defined above and
$(X_{\delta_N},X_{2\delta_N},\ldots,X_{N\, \delta_N})$ is known,
with $N (\delta_N)^2 \limiteNN 0$ and $N \, \delta_N \limiteNN
\infty$. Under Assumptions $W_3$ and $C$, using $\tilde
S^{k'}_{k}(a_N)$ instead of $S^{k'}_{k}(a_N)$, if $ \frac \mu
\lambda < \frac {f_{max}}{f_{min}}$ and $r_i=\frac
{f_{min}}\lambda+\frac i \ell \big (\frac {f_{max}}\mu-\frac
{f_{min}}\lambda\big )$ for $i=1,\ldots,\ell$ with $a_N=1$ for all
$N \in \Inte$, then limit theorems (\ref{TLCS}), (\ref{lim_epsi})
and (\ref{TLClogS}) hold with $\alpha_j^*=2H_j^*+1$ and
$\beta_j^*=\log \big (-\frac {\sigma_j^{*2}} 2 \int_\Real \big
|\widehat \psi(u)\big |^2\, |u|^{-1-2H_j^*} du\big)$ for all
$j=0,1,\ldots,m$, for all $(p,q)\in \{1,\ldots,\ell\}$,
\begin{eqnarray}  \label{TLCdis}
\gamma^{(j)}_{pq}=\frac {2}{(1-2w)\left(r_p\, r_q \right)^{2H_j^*}}
\int _{\Real}\Big( \frac { \int _{\Real}\overline{\widehat{\psi}}(
r_p\, \xi) \widehat{\psi}( r_q \, \xi) \,
|\xi|^{-1-2H_j^*}e^{-iu\xi}d \xi } {\int_\Real \big |\widehat
\psi(u)\big |^2\, |u|^{-1-2H_j^*} du} \Big )^2du.
\end{eqnarray}
\end{propo}
Theorem \ref{convprob} can be applied to a piecewise locally
fractional Gaussian process without conditions on parameters
$H_j^*$. Thus,
\begin{cor}\label{cor_local}
Under assumptions of Proposition \ref{Theo:TCL:discret} and
Assumption C, then for all $0<\kappa<\frac 1 2$, if
$\delta_N=N^{-1/2-\kappa}$ and $v_N=N^{1/2-\kappa}$ then
(\ref{conv_tau}), (\ref{conv_theta}),
 (\ref{estim2}) and (\ref{test}) hold.
\end{cor}
Therefore the convergence rate of $\underline{\widehat \tau}$ to
$\underline{\tau}^*$ (in probability) is as well close to $N^{1/2}$
as one wants. Estimators $\tilde H_j$ and ${\overline H_j}$ converge
to the parameters $H_j^*$ following a central limit theorem with a
rate of convergence $N^{1/4-\kappa/2}$ for $0<\kappa$ as small as
one wants.
\subsection{Application to heart rate's time series}
\begin{center}
\begin{figure}[h]
\begin{center}
\includegraphics[width=10 cm,height=4 cm]{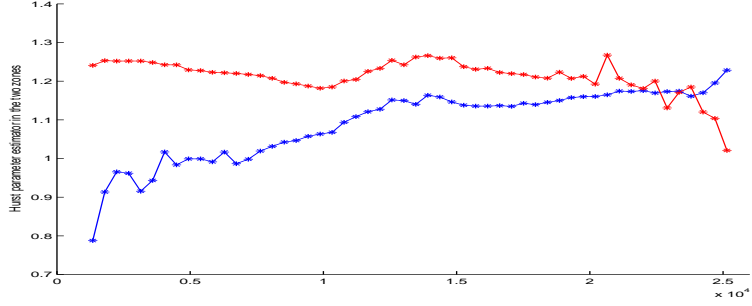}
\end{center}
\caption{Evolution of local-fractality parameter estimators
(observed for HR series of one athlete) in the two zones when the
change point varies in time.} \label{evolhurst}
\end{figure}
\end{center}

The study of the regularity of physiological data and in particular
the heartbeat signals have received much attention by several
authors (see for instance \cite{pengh}, \cite{pengm} or
\cite{absil}). These authors studied HR series for healthy subjects
and subjects with heart disease. In \cite{bardet2}, a piecewise
locally fractional Brownian motion is studied for modeling the
cumulative HR data during three typical phases (estimated from
Lavielle's algorithm) of the race (beginning, middle and end). The
local-fractality parameters are estimated with wavelet analysis. The
conclusions obtained are relatively close to those obtained by Peng.
\emph{et al.}. Indeed we remarked that the local-fractality
parameter increases thought the race phases which may be explained
with fatigue appearing during the last phase of the marathon. In
this paper, we tray to unveil in which instants the behaviour of HR
data changes. The following Table \ref{Tableath} presents the
results for the detection of one change point.
\begin{table}[h]
\renewcommand{\arraystretch}{1.2}
\begin{center}
\begin{tabular}{cccccc}
 $ $ &$\widehat \tau_1$ &
$\tilde H_0$ & $\tilde H_1$ & $T^{(0)}$ & $T^{(1)}$\\
\hline \hline
Ath1 & 0.0510 & 0.7880 & 1.2376  & 1.0184 & 1.0562\\
Ath2 & 0.4430 & 1.3470 & 1.4368  & 5.0644 & 1.5268\\
Ath3 & 0.6697 & 0.9542 & 1.2182  & 0.7836 & 0.9948\\
Ath4 & 0.4856 & 1.1883 & 1.2200  & 2.8966 & 1.2774\\
Ath5 & 0.8715 & 1.1512 & 1.3014  & 0.7838 & 0.8748\\
Ath6 & 0.5738 & 1.1333 & 1.1941  & 2.2042 & 0.7464\\
Ath7 & 0.3423 & 1.1905 & 1.1829  & 0.4120 & 1.5598\\
Ath8 & 0.8476 & 1.0222 & 1.2663  & 3.1704 & 0.5150\\
Ath9 & 0.7631 & 1.4388 & 1.3845  & 9.6574 & 0.5714\\
\hline \hline
\end{tabular}
\end{center}
\caption{Estimated change points $\tau_1$, parameters $H_0$, $H_1$
and goodness-of-fit test statistics ($T^{(0)}$ for the first zone
and $T^{(1)}$ for the second one) in the case of one change point
observed in HR series of different athletes.} \label{Tableath}
\end{table}
\newline It is noticed that the estimator of the local-fractality
parameter is generally larger on the second zone than on the first
although the detected change point differs from an athlete to
another (only the case of Athlete 1 seems not to be relevant). This
result is very interesting and confirms our conclusions in
\cite{bardet2}. Whatever is the position of change point, the
estimation of the local-fractality parameter is larger in the second
segment than in the first segment (see the example of HR data
recorded for one athlete in Figure \ref{evolhurst}).
In general, the goodness-of-fit tests, with values $T^{(0)}$ and
$T^{(1)}$, are less than $\chi^2_{95\%}(4)=9.4877$ (except $T^{(0)}$
for Ath9) when $\ell=6$. So, the HR data trajectory in the both
zones seems to be correctly modeled with a stationary locally
fractional Gaussian trajectory.
\section{Proofs} \label{proofs}
Before establishing the proof of Theorem \ref{convprob} an important
lemma can be stated:
\begin{lemma}\label{lemme1}
Let $k \in \Inte \setminus \{0,1\}$, $(\gamma_i)_{1\leq i \leq k}
\in (0,\infty)^k$ and $\alpha_1>\alpha_2>\cdots>\alpha_k$ be $k$
ordered real numbers. For $(\alpha,\beta) \in \Real^2$, consider the
function $f_{\alpha,\beta}:x\in \Real \mapsto \Real$ such that
\[
f_{\alpha,\beta}(x):=\alpha \, x+\beta-\log\Big(\sum_{q=1}^{k}
\gamma_q \exp \big(\alpha_q\, x \big)\Big) ~~\mbox{for}~ x \in
\Real.
\]
Let $0<t_1<\cdots<t_\ell$ with $\ell \in \Inte \setminus \{0,1,2\}$
and $(u_n)_{n \in \Inte}$ be a sequence of real numbers such that
there exists $m \in \Real$ satisfying $u_n\geq m$ for all $n \in
\Inte$. Then there exists $C>0$ not depending on $n$ such that
$$
\inf_{(\alpha,\beta) \in \Real^2} \sum_{i=1}^\ell \big |
f_{\alpha,\beta} \big (\log(u_n)+t_i\big ) \big |^2 \geq C \,
\min\big (1,|u_n|^{2(\alpha_2-\alpha_1)}\big).
$$
\end{lemma}
~\\
\noindent {\em \underline{Proof of Lemma \ref{lemme1}}:} For all
$(\alpha,\beta) \in \Real^2$, the function $f_{\alpha,\beta} $ is a
$\mathcal{C}^\infty(\Real)$ function and
\[
\frac{\partial^2 }{\partial
x^2}f_{\alpha,\beta}(x)=-\frac{\sum_{q=1}^{k-1}
\gamma_q\gamma_{q+1}(\alpha_q-\alpha_{q+1})^2\exp
\big((\alpha_q+\alpha_{q+1})\, x\big)}{\Big(\sum_{q=1}^{k}
\gamma_q\exp \big(\alpha_q \, x\big)\Big)^2}<0.
\]
Therefore the function $f_{\alpha,\beta}$ is a concave function such
that $\sup_{(\alpha,\beta) \in \Real^2}\frac{\partial^2 }{\partial
x^2}f_{\alpha,\beta}(x)<0$ (not depending on $\alpha$ and $\beta$)
and for all $(\alpha,\beta) \in \Real^2$, $f_{\alpha,\beta}$
vanishes in $2$ points at most. Thus, since $\ell \geq 3$ and $\big
(x+t_i\big )_i$ are distinct points, for all $x \in \Real$, it
exists $C(x)>0$ not depending on $\alpha$ and $\beta$ such that
$$
\inf_{(\alpha,\beta) \in \Real^2} \sum_{i=1}^\ell \big |
f_{\alpha,\beta} \big (x+t_i\big ) \big |^2 \geq C(x).
$$
Therefore, since for all $M\geq 0$,
\begin{eqnarray}\label{ineg_borne}
\inf_{x \in [-M,M]} \left \{ \inf_{(\alpha,\beta) \in \Real^2}
\sum_{i=1}^\ell \big | f_{\alpha,\beta} \big (x+t_i\big ) \big |^2
\right \}  \geq \inf_{x \in [-M,M]} \big \{ C(x) \big \}>0.
\end{eqnarray}
Moreover, if $u_n \to +\infty$,
\begin{eqnarray*}
\log\Big(\sum_{q=1}^{k} \gamma_q \exp \big(\alpha_q\, \log(u_n)
\big)\Big)&=& \log\Big( \gamma_1 \exp \big(\alpha_1\, \log(u_n)
\big)+\gamma_2 \exp \big(\alpha_2\, \log(u_n) \big)\big (1 +o(1)\big
)
\Big)\\
&=& \log (\gamma_1)+ \alpha_1 \log(u_n)+\gamma_2 \exp
\big((\alpha_2-\alpha_1)\, \log(u_n) \big)\big (1 +o(1)\big ).
\end{eqnarray*}
Thus, for $n$ large enough, \begin{eqnarray}\label{ineg_f} \frac 1 2
\, \gamma_2 \, u_n^{\alpha_2-\alpha_1} \leq  \Big |
\log\Big(\sum_{q=1}^{k} \gamma_q \exp \big(\alpha_q\, \log(u_n)
\big)\Big)- \log (\gamma_1)+ \alpha_1 \log(u_n)\Big | \leq 2 \,
\gamma_2 \, u_n^{\alpha_2-\alpha_1}.
\end{eqnarray}
Therefore, for all $(\alpha,\beta) \in \Real^2$,
\begin{multline*}
\big | f_{\alpha,\beta} \big (\log(u_n)+t_i\big ) \big |^2 =\big |
f_{\alpha_1,\log ( \gamma_1)} \big (\log(u_n)+t_i\big )\big | ^2 +
\Big |(\log (\gamma_1)-\beta)+ (\alpha_1-\alpha)
(\log(u_n)+t_i)\Big |^2 \\
 - 2f_{\alpha_1,\log ( \gamma_1)} \big
(\log(u_n)+t_i\big ) \times \Big ((\log (\gamma_1)-\beta)+
(\alpha_1-\alpha) (\log(u_n)+t_i)\Big ).
\end{multline*}
Using inequalities (\ref{ineg_f}), $\displaystyle \frac 1 4 \,
\gamma_2^2 \, u_n^{2(\alpha_2-\alpha_1)} \leq \big |
f_{\alpha_1,\log ( \gamma_1)} \big (\log(u_n)+t_i\big )\big | ^2
\leq 4 \, \gamma_2^2 \, u_n^{2(\alpha_2-\alpha_1)} $ and for all
$(\alpha,\beta) \in \Real^2$,  $\displaystyle \lim_{n \to \infty}
f_{\alpha_1,\log ( \gamma_1)} \big (\log(u_n)+t_i\big ) \times \Big
((\log (\gamma_1)-\beta)+ (\alpha_1-\alpha) (\log(u_n)+t_i)\Big
)=0$. Then, for all $\displaystyle (\alpha,\beta) \neq
(\alpha_1,\log ( \gamma_1))$, $\displaystyle \lim_{n \to \infty}
\big | f_{\alpha,\beta} \big (\log(u_n)+t_i\big ) \big |^2=\infty$.
Consequently,  for $n$ large enough,
\begin{eqnarray*}
\inf_{(\alpha,\beta) \in \Real^2} \sum_{i=1}^\ell \big |
f_{\alpha,\beta} \big (\log(u_n)+t_i\big ) \big |^2  &\geq &\frac 1
2 \, \sum_{i=1}^\ell \big | f_{\alpha_1,\log ( \gamma_1)} \big
(\log(u_n)+t_i\big ) \big |^2 \\
 &\geq &\frac
1 8 \, \gamma_2^2 \,\sum_{i=1}^\ell
(u_n+t_i)^{2(\alpha_2-\alpha_1)}\\
&\geq & C\, u_n^{2(\alpha_2-\alpha_1)},
\end{eqnarray*}
which combined with (\ref{ineg_borne}) achieves the proof.
\hfill $\Box$ \\
\\
\noindent {\em \underline{Proof of Theorem \ref{convprob}}:} Let
$\displaystyle w_N=\frac {N\delta_N} {v_N}$, $k_j^*=[N\delta_N
\tau_j^*]$ for $j=1,\ldots,m$ and
$$V_{\eta \, w_N}=\{(k_j)_{1\leq j\leq
m},~\max_{j\in{1,\ldots,m}}|k_j-k^*_j|\geq \eta \,w_N\}.
$$
Then, for $N\delta_N$ large enough,
\begin{eqnarray}
\nonumber \Prob \Big (\frac {N\delta_N} {w_N} \,
\|\underline{\tau}^*-\underline{\widehat \tau}\|_m \geq \eta \Big )&
\simeq & \Prob(\max_{j\in{1,\ldots,m}}|\widehat k_j-k^*_j|\geq \eta
\,w_N)\\
\nonumber &= &\Prob \Big(\min_{(k_j)_{1\leq j\leq m}\in V_{\eta
\,w_N}} G_N\big((k_j)_{1\leq j\leq m}\big)\leq\min_{(k_j)_{1\leq
j\leq m}\notin V_{\eta \,w_N}} G_N\big((k_j)_{1\leq j\leq
m}\big)\Big)\\
\label{prob} &\leq &\Prob \Big(\min_{(k_j)_{1\leq j\leq m}\in
V_{\eta \,w_N}} G_N\big((k_j)_{1\leq j\leq m}\big)\leq G_N\big (
(k_j^*)_{1\leq j\leq m}\big)\Big).
\end{eqnarray}
For $j=\{0,\ldots,m\}$ and $0=k_0<k_1<\ldots<k_m<k_{m+1}=N\delta_N$,
let
\begin{itemize}
\item $Y_{k_j}^{k_{j+1}}:=\big(\log\big(S^{k_{j+1}}_{k_j}(r_i\cdot
a_N)\big)\big)_{1\leq i\leq \ell}$,

\item $\Theta^{k_{j+1}}_{k_j}=\left(
\begin{array}{c}
\alpha_j \\
\log \beta_j\\
\end{array} \right)$,~~~ $\widehat \Theta^{k_{j+1}}_{k_j}=\left(
\begin{array}{c}
\widehat \alpha_j \\
\log \widehat\beta_j\\
\end{array} \right)$~~ and~~$\Theta^*_j= \left(
\begin{array}{c}
\alpha^*_j \\
\log \beta^*_j\\
\end{array}
\right)$.
\end{itemize}
1/ Using these notations, $\displaystyle{G_N\big((k_j)_{1\leq j\leq
m}\big)=\sum_{j=0}^m \parallel Y^{k_{j+1}}_{k_j}-L_{a_N} \cdot
\widehat\Theta^{k_{j+1}}_{k_j}\parallel^2}$, where $\parallel \cdot
\parallel$ denotes the usual Euclidean norm in $\Real ^\ell$. Then,
with $I_\ell$ the $(\ell\times \ell)$-identity matrix
\begin{eqnarray*}
G_N\big ( (k_j^*)_{1\leq j\leq m}\big)&=&\sum_{j=0}^m
\parallel Y^{k^*_{j+1}}_{k^*_j}-L_{a_N} \cdot \Theta^*_j\parallel^2\\
&=&\sum_{j=0}^m \Big\| (I_\ell-P_{L_{a_N}})\cdot
Y_{k_j^*}^{k_{j+1}^*}\Big\|^2 \mbox{~~~~~~~~~~~~~with $P_{L_{a_N}}=L_{a_N}\cdot(L_{a_N}'\cdot L_{a_N})^{-1}\cdot L_{a_N}'$}\\
&=&  \sum_{j=0}^m \frac {a_N }{k^*_{j+1}-k^*_j} \Big\|
(I_\ell-P_{L_{a_N}})\cdot \big (\varepsilon_i^{(N)}(k^*_j,k^*_{j+1})
\big )_{1\leq i\leq \ell} \Big\|^2~~~~~~~~~~~~~~\mbox{from (\ref{TLClogS})}\\
& \leq &\frac 1 {\min_{0\leq j \leq m}(\tau^*_{j+1}-\tau^*_j)} \cdot
\frac {a_N }{N \, \delta_N} \, \sum_{j=0}^m  \Big\|  \big
(\varepsilon_i^{(N)}(k^*_j,k^*_{j+1}) \big )_{1\leq i\leq
\ell}\Big\|^2.
\end{eqnarray*}
Now, using the limit theorem (\ref{lim_epsi}), $\displaystyle \Big\|
\big (\varepsilon_i^{(N)}(k^*_j,k^*_{j+1}) \big )_{1\leq i\leq
\ell}\Big\|^2 \limitedistl \Big \|
\mathcal{N}(0,\Gamma(r_1,\ldots,r_\ell))\Big\| ^2$ since
$k^*_{j+1}-k^*_j \sim N\delta_N(\tau_{i+1}^*-\tau_i^*) \limiteNN
\infty$, and thus
\begin{eqnarray}\label{G*}
G_N\big ( (k_j^*)_{1\leq j\leq m}\big)=O_P\Big (\frac {a_N }{N \,
\delta_N} \Big ),
\end{eqnarray}
where $\xi_N=O_P(\psi_N)$ as $N \to \infty$ is written, if for all
$\rho>0$, there exists $c>0$, such as $ P\Big (|\xi_N| \leq c \cdot
\psi_N \Big ) \geq 1 -\rho$ for all sufficiently large $N$.\\
~\\
2/ Now, set $(k_j)_{1\leq j\leq m}\in V_{\eta \,w_N}$. Therefore,
for $N$ and $N\,\delta_N$ large enough, there exists $j_0 \in
\{1,\ldots,m\}$ and $(j_1,j_2) \in \{1,\ldots,m\}^2$ with $j_1 \leq
j_2$ such that $k_{j_0}\leq k^*_{j_1}-\eta \,w_N$ and $k_{j_0+1}\geq
k^*_{j_2}+\eta \,w_N$. Thus,
\[
G_N\big((k_j)_{1\leq j\leq m}\big )\geq \big \|
Y_{k_{j_0}}^{k_{{j_0}+1}}-L_{a_N}
\widehat\Theta_{k_{j_0}}^{k_{{j_0}+1}}\big \| ^2.
\]
Let $\Omega^*:= \big (\Omega^*_i\big )_{1\leq i \leq
\ell}$ be the vector such that\\
\begin{eqnarray*}
\Omega^*_i&:=&{\frac{k^*_{j_1}-k_{j_0}}{k_{{j_0}+1}-k_{j_0}}}\,
\beta^*_{j_1-1} \exp \big(\alpha^*_{j_1-1}\log (r_i \, a_N)
\big)+\sum_{j=j_1}^{j_2-1}
{\frac{k^*_{j+1}-k^*_{j_1}}{k_{{j_0}+1}-k_{j_0}}}\beta^*_{j} \exp
\big(\alpha^*_{j}\log (r_i \, a_N) \big)\\
&& \hspace{7cm}
+{\frac{k_{j_0+1}-k^*_{j_2}}{k_{{j_0}+1}-k_{j_0}}}\beta^*_{j_2} \exp
\big(\alpha^*_{j_2}\log (r_i \, a_N) \big).
\end{eqnarray*}
Then, \begin{eqnarray}\label{majG}G_N\big((k_j)_{1\leq j\leq m}\big
) \geq
\parallel
Y_{k_{j_0}}^{k_{{j_0}+1}}-\big (\log \Omega^*_i\big )_{1\leq i \leq
\ell}\parallel^2+\parallel \big (\log \Omega^*_i\big )_{1\leq i \leq
\ell}-L_{a_N} \cdot
\widehat\Theta_{k_{j_0}}^{k_{{j_0}+1}}\parallel^2 +2 \, Q,
\end{eqnarray}
with $Q=\big (Y_{k_{j_0}}^{k_{{j_0}+1}}-\big (\log \Omega^*_i\big
)_{1\leq i \leq \ell}\big )' \cdot \big ( \big (\log \Omega^*_i\big
)_{1\leq i \leq \ell}-L_{a_N} \cdot
\widehat\Theta_{k_{j_0}}^{k_{{j_0}+1}}\big ).$ \\
In the one hand, with $S^{k'}_{k}(\cdot)$ defined in (\ref{Sk}),
\begin{eqnarray*}
S_{k_{j_0}}^{k_{{j_0}+1}}(r_i \,
a_N)={\frac{k^*_{j_1}-k_{j_0}}{k_{{j_0}+1}-k_{j_0}}}\,
S_{k_{j_0}}^{k^*_{j_1}}(r_i \, a_N)+\sum_{j=j_1}^{j_2-1}
{\frac{k^*_{j+1}-k^*_{j_1}}{k_{{j_0}+1}-k_{j_0}}}S_{k^*_{j}}^{k^*_{j+1}}(r_i
\, a_N)+{\frac{k_{j_0+1}-k^*_{j_2}}
{k_{{j_0}+1}-k_{j_0}}}S_{k_{j_2}}^{k^*_{j_0+1}}(r_i \, a_N).
\end{eqnarray*}
Using the central limit theorems (\ref{TLClogS}), for $N$ and
$N\,\delta_N$ large enough,
\begin{eqnarray*}
\Esp \Big [\big (S_{k_{j_0}}^{k_{{j_0}+1}}(r_i \, a_N)-\Omega^*_i
\big )^2 \Big ] &\leq & m \, \left ( \Big (
{\frac{k^*_{j_1}-k_{j_0}}{k_{{j_0}+1}-k_{j_0}}}\Big )^2 \,\Esp \Big
[\big ( S_{k_{j_0}}^{k^*_{j_1}}(r_i \, a_N) - \beta^*_{j_1-1}
\big(r_i \, a_N \big)^{\alpha^*_{j_1-1}}\big )^2 \Big ] \right . \\
&& +\sum_{j=j_1}^{j_2-1} \Big (
{\frac{k^*_{j+1}-k^*_{j_1}}{k_{{j_0}+1}-k_{j_0}}}\Big ) ^2  \,\Esp
\Big [\big ( S_{k^*_{j}}^{k^*_{j+1}}(r_i \, a_N)-\beta^*_{j}
\big(r_i \, a_N \big)^{\alpha^*_{j}}\big )^2 \Big ] \\
&& \left . +\Big ( {\frac{k_{j_0+1}-k^*_{j_2}}
{k_{{j_0}+1}-k_{j_0}}}\Big ) ^2 \,\Esp \Big [\big
(S_{k_{j_2}}^{k^*_{j_0+1}}(r_i \, a_N)-\beta^*_{j_2} \big(r_i \,
a_N \big)^{\alpha^*_{j_2}}\big )^2 \Big ]\right )\\
\Longrightarrow ~~~\Esp \Big [\big (\frac
{S_{k_{j_0}}^{k_{{j_0}+1}}(r_i \, a_N)}{\Omega^*_i }-1\big )^2 \Big
] &\leq & \frac {m\,\gamma^2}  {\Omega^*_i } \, {a_N} \,
 \, \Big (\frac 1 {k^*_{j_1}-k_{j_0}}+\sum_{j=j_1}^{j_2-1}
\frac 1 {k^*_{j+1}-k^*_{j_1}}+\frac 1 {k^*_{j_0+1}-k_{j^*_2}} \Big
)\\
&\leq & C \,\frac {a_N} {\eta\,  w_N },
\end{eqnarray*}
with $\gamma^2=\max_{i,j}\{\gamma_{ii}^{(j)}\}$ (where
$(\gamma_{pq}^{(j)})$ is the asymptotic covariance  of vector $
\varepsilon_p^{(N)}(k,k')$ and $ \varepsilon_q^{(N)}(k,k')$) and
$C>0$ not depending on $N$. Therefore, for $N$ large enough, for all
$i=1,\ldots,\ell$,
$$
\Esp \Big [\big (\log (S_{k_{j_0}}^{k_{{j_0}+1}}(r_i \,
a_N))-\log(\Omega^*_i)\big )^2 \Big ] \leq  2 \, C \,\frac {a_N}
{\eta\,  w_N}.
$$
Then we deduce with Markov Inequality that
\begin{eqnarray}\label{Q1}
\parallel
Y_{k_{j_0}}^{k_{{j_0}+1}}-\big (\log \Omega^*_i\big )_{1\leq i \leq
\ell}\parallel^2=O_P\Big(\frac {a_N} {\eta\,  w_N}\Big).
\end{eqnarray}
>From the other hand,
\begin{eqnarray*}
\parallel \big (\log \Omega^*_i\big )_{1\leq i
\leq \ell}-L_{a_N} \cdot
\widehat\Theta_{k_{j_0}}^{k_{{j_0}+1}}\parallel^2=
\sum_{i=1}^{\ell}\Big(\big(\widehat\alpha_{j_0}\log (r_i \, a_N)
+\log\widehat \beta_{j_0}\big )-\log \Omega^*_i \Big)^2.
\end{eqnarray*}
Define $\displaystyle
\gamma_{1}:={\frac{k^*_{j_1-1}-k_{j_0}}{k_{{j_0}+1}-k_{j_0}}}\cdot
\beta^*_{j_1-1}$, for all $p\in \{0,1,\ldots,j_2-j_1-1\}$,
$\displaystyle
\gamma_{p}:={\frac{k^*_{j_1+p}-k^*_{j_1+p-1}}{k_{{j_0}+1}-k_{j_0}}}\cdot
\beta^*_{j_1+p-1}$ and $\displaystyle
\gamma_{j_2-j_1+1}:={\frac{k_{j_0+1}-k^*_{j_2}}{k_{{j_0}+1}-k_{j_0}}}\cdot
\beta^*_{j_2}$. Then, using Lemma \ref{lemme1}, one obtains
$$
\inf_{\alpha,\beta}\Big \{ \sum_{i=1}^{\ell}\Big(\big(\alpha\log (
r_i\, a_N)+\log \beta \big )-\log \Omega^*_i \Big)^2\Big \} \geq C
\, \min\big (1_,,\, |a_N|^{2(\alpha^*_{(2)}-\alpha^*_{(1)})}\big),
$$
where $C>0$ and $\alpha^*_{(1)}=\max_{j=j_1-1,\ldots,j_2}
\alpha_j^*$, $\alpha^*_{(2)}=\max_{j=j_1-1,\ldots,j_2,~j\neq (1)}
\alpha_j^*$. As a consequence, for satisfying all possible cases of
$j_0$, $j_1$ and $j_2$, one obtains
\begin{eqnarray}\label{Q2}
\parallel \big (\log \Omega^*_i\big )_{1\leq i
\leq \ell}-L_{a_N} \cdot
\widehat\Theta_{k_{j_0}}^{k_{{j_0}+1}}\parallel^2 \geq C \,
|a_N|^{2(\min_i \alpha^*_{i}-\max_i \alpha^*_{i})}.
\end{eqnarray}
Finally, using Cauchy-Schwarz Inequality,
\begin{eqnarray*}
Q&\leq &\Big ( \big \| Y_{k_{j_0}}^{k_{{j_0}+1}}-\big (\log
\Omega^*_i\big )_{1\leq i \leq \ell}\big \|^2 \cdot \big \| \big
(\log \Omega^*_i\big )_{1\leq i \leq \ell}-L_{a_N} \cdot
\widehat\Theta_{k_{j_0}}^{k_{{j_0}+1}}\big |^2 \Big )^{1/2}
\end{eqnarray*}
Therefore, using (\ref{Q1}) and (\ref{Q2}), since under assumptions
of Theorem \ref{convprob},
$$\displaystyle \frac {a_N} {\eta\, w_N}=o\Big ( |a_N|^{2(\min_i \alpha^*_{i}-\max_i
\alpha^*_{i})}\Big ),$$ then
\begin{eqnarray}\label{Q}
Q=o_P\Big(\parallel \big (\log \Omega^*_i\big )_{1\leq i \leq
\ell}-L_{a_N} \cdot
\widehat\Theta_{k_{j_0}}^{k_{{j_0}+1}}\parallel^2\Big).
\end{eqnarray}
We deduce from relations (\ref{majG}), (\ref{Q1}), (\ref{Q2}) and
(\ref{Q}) that
$$\Prob \Big (\min_{(k_j)_{1\leq j\leq m}\in V_{\eta \,w_N}} G_N\big((k_j)_{1\leq j\leq
m}\big )\geq \frac C 2 \, |a_N|^{2(\min_i \alpha^*_{i}-\max_i
\alpha^*_{i})}\Big ) \limiteNN 1.$$ $\hfill \Box$\\
\\
\noindent {\em \underline{Proof of Theorem \ref{convparam}}:} From
Theorem \ref{convprob}, it is clear that
\begin{eqnarray}\label{inclus}
\Prob\big([\tilde k_j,\tilde
k_{j}']\subset[k_j^*,k^*_{j+1}]\big)\limiteNN 1~~\mbox{and}~~ \frac
{\tilde k_{j}'-\tilde
k_j}{N\delta_N(\tau_{j+1}^*-\tau_j^*)}\limiteprob 1.
\end{eqnarray}
Now, for $j=0,\ldots,m$, $(x_i)_{1\leq i\leq\ell}\in\Real^\ell$ and
$0<\varepsilon <1$, let $A_j$ and $B_j$ be the events such that
\begin{multline*}
A_j:=\Big\{[\tilde k_j,\tilde
k_{j}']\subset[k_j^*,k^*_{j+1}]\Big\}\bigcap \left \{ \Big |\frac
{\tilde k_{j}'-\tilde k_j}{N\delta_N(\tau_{j+1}^*-\tau_j^*)}-1 \Big|
\leq \varepsilon \right \} \\
~~\mbox{and}~~B_j:=\left \{\sqrt {\frac {\tilde k_{j}'-\tilde
k_j}{a_N}} \Big(Y_{\tilde k_j}^{\tilde k_{j}'}-L_{a_N} \cdot
\Theta^*_j\Big) \in \prod_{i=1}^{\ell}(-\infty,x_i]\right \}
\end{multline*}
First, it is obvious that
\begin{eqnarray}\label{AB}
\Prob(A_j)\Prob(B_j\mid A_j)\leq\Prob(B_j)\leq\Prob(B_j\mid
A_j)+1-\Prob(A_j).
\end{eqnarray}
Moreover, from (\ref{TLCS}),
\begin{eqnarray*}
\Prob(B_j\mid A_j)&=&\Prob\Big( \big(\varepsilon_i^{(N)}(\tilde
k_j,\tilde k_{j}')\big)_{1\leq
i\leq\ell}\in\prod_{i=1}^{\ell}(-\infty,x_i]\mid A_j\Big)\\
&\limiteNN &
\Prob\Big(\Nor\big(0,\Gamma^{(j)}(\alpha_j^*,r_1,\ldots,r_\ell)\big)
\in\prod_{i=1}^{\ell}(-\infty,x_i] \Big ).
\end{eqnarray*}
Using (\ref{inclus}), it is straightforward that $\Prob(A_j)
\limiteNN 1$. Consequently,
$$
\Prob(B_j)\limiteNN
\Prob\Big(\Nor\big(0,\Gamma^{(j)}(\alpha_j^*,r_1,\ldots,r_\ell)\big)
\in\prod_{i=1}^{\ell}(-\infty,x_i] \Big )
$$
and therefore $\displaystyle \sqrt {\frac {\tilde k_{j}'-\tilde
k_j}{a_N}} \Big(Y_{\tilde k_j}^{\tilde k_{j}'}-L_{a_N} \cdot
\Theta^*_j\Big) \limitedistl
\Nor\big(0,\Gamma^{(j)}(\alpha_j^*,r_1,\ldots,r_\ell)\big )$. Now
using again (\ref{inclus}) and Slutsky's Lemma one deduces
$$ \sqrt
{\frac {\delta_N \, \big(N(\tau_{j+1}^*-\tau_j^*)\big)}{a_N}}
\Big(Y_{\tilde k_j}^{\tilde k_{j}'}-L_{a_N} \cdot \Theta^*_j\Big)
\limitedistl
\Nor\big(0,\Gamma^{(j)}(\alpha_j^*,r_1,\ldots,r_\ell)\big ).
$$
Using the expression of $\tilde\Theta_j$ as a linear application of
$Y_{\tilde k_j}^{\tilde k_{j}'}$, this achieves the proof of Theorem
\ref{convparam}. $\hfill \Box$

\bibliographystyle{amsplain}

\end{document}